\documentclass[ 10pt, a4paper,leqno]{amsart}
\usepackage[francais,english]{babel}
\theoremstyle{plain}
 \def\dem{\begin{proof}[Preuve]} 
  \def\CQFD{\end{proof}} 
    \def\N{ \mbox{I\hspace{-.15em}N}}
       \def\Z{ \mbox{Z\hspace{-.3em}Z}}   
        
         \def\P{ \mbox{I\hspace{-.17em}P}}

               \newtheorem {thm}{Th\'eor\`eme\hspace{2pt}} 
               
              \newtheorem {pro}{Proposition\hspace{2pt} }  
                
           \newtheorem {cor}{Corollaire\hspace{2pt}}  
           \newtheorem {defi}{Definition\hspace{2pt}}  
         
        \newtheorem {lem}{Lemme\hspace{2pt}}
\author[Fabre]{Bruno Fabre}

\address{22, rue Emile Dubois \\
 75014 Paris\\
 France}
 \email{bruno.fabre9@wanadoo.fr}
\title{Fonctions de Hilbert et g\'eom\'etrie \\ Hilbert functions and geometry}
 
\begin{document}
\begin{abstract}
This note is devoted to the study of the links between
 the Hilbert function of a subscheme $X$ of the projective space $\P_n$ and its geometric properties. 
 We will assume $X$ to be arithmetically Cohen-Macaulay (ACM).
  This allows us to characterize the Hilbert function $\phi_X$ of $X$ by an increasing sequence of $d$ integers $(m_0,\dots,m_{d-1})$,
   called the (absolute) {\it characteristic sequence} of $X$, $d$ being the degree of $X$. 
  If $Y$ is an ACM hypersurface of $X$, we characterize the Hilbert function of $X$
   by a increasing sequence of $d$ integers, called the {\it relative characteristic sequence}
  of $Y$ in $X$. We study properties of these sequences, and study
   in this context, on a Gorenstein curve $X$, linear systems with maximal dimension with
   respect to their degree.
   \end{abstract}

\subjclass{14C20, 13A02}
\keywords{Fonctions de Hilbert, groupes de points, syst\`emes lin\'eaires. }
\thanks{ Je remercie C. Ciliberto pour son accueil \`a l'universit\'e Tor Vergata, 
 ainsi que les discussions qu'il a bien voulu m'accorder.
   Je remercie \'egalement d'autres math\'ematiciens pour des discussions fructueuses, dont P. Mazet, C. Peskine, J.-L. Sauvageot, M. Chardin. }
\maketitle

\frenchspacing
Cette note est consacr\'ee \`a l'\'etude des relations entre la fonction de Hilbert $\phi_{X}$
 d'un sous-sch\'ema projectif $X\subset\P_{n}$ (i.e. de son c\^one associ\'e) et ses propri\'et\'es g\'eom\'etriques.
  Lorsque $X$ est un groupe de points du plan, l'\'etude est largement avanc\'ee, notamment gr\^ace \`a l'introduction du {\it carac\-t\`ere num\'erique}
   de Gruson-Peskine (\cite{Peskine}). On a montr\'e en particulier dans \cite{Fabre} comment on pouvait retrouver
 la description g\'eom\'etrique des syst\`emes lin\'eaires de dimension maximale sur une courbe plane $X$,
  donn\'ee par Ciliberto dans \cite{Ci2} lorsque $X$ est lisse, gr\^ace \`a ce caract\`ere num\'erique.
   Pour $n\ge 3$, le probl\`eme \'etait pos\'e depuis longtemps de trouver une g\'en\'eralisa\-tion ad\'equate du 
   caract\`ere num\'erique pour les groupes de points dans $\P_n$, en particulier pour l'\'etude des groupes de points et des syst\`emes lin\'eaires 
   sur les courbes alg\'ebriques de $\P_n$. On propose ici une telle g\'en\'eralisation, 
   en introduisant le concept de suite caract\'eristique relative d'un groupe de points sur une courbe alg\'ebrique arithm\'etiquement Cohen-Macaulay (ACM). 
   La motivation originelle \'etait de g\'en\'eraliser cette description g\'eom\'etrique des syst\`emes lin\'eaires de dimension maximale obtenue pour les courbes planes, 
  aux courbes de Gorenstein. On pense ici en particulier aux th\'eor\`emes obtenus
    pour les intersections compl\`etes par C. Ciliberto et R. Lazarsfeld dans \cite{Ci1} et par B. Basili dans \cite{Basili}.
L'objectif n'est pas encore atteint, mais on a les r\'esultats suivants.

 Etant donn\'e un sous-sch\'ema projectif ACM
  $X\subset \P_n$ non-d\'eg\'en\'er\'e de dimension $m$ et de degr\'e $d$, 
  l'anneau gradu\'e projetant $A_X:=A_{n}/I_{X}$ de $X$ s'\'ecrit $A_X=\sum_{i=0}^{d-1}{A_m[-m_i]}$ ($A_{m}=k[X_{0},\dots,X_{m}¥]¥$¥), 
  les $m_i$ \'etant une suite croissante de $d$ entiers caract\'erisant la fonction de Hilbert de $X$, appel\'ee sa {\it suite caract\'eristique}. On pose aussi:
  $l_{i}:=card\{j/m_{j}=i\}$. On a :
 
 1. $l_0=1$, $l_{1}=p:=n-m$. Si $l_i=1$, alors $l_{i+1}\ge p$. De plus, $l_{i}$ v\'erifie la condition de croissance de Macaulay (cf. \cite{Gotzmann}¥);
  en particulier si $l_{i}=0, l_{i+1}=0$.
 
 2. Si $X$ est contenu dans $X'$, la suite caract\'eristique de $X$ est contenue dans celle de $X'$.
 
 Soit maintenant $Y\subset X$ une hypersurface ACM de $X$. On peut caract\'eriser la fonction de Hilbert
  par une suite croissante de $d$ entiers $(n_i)$, appel\'ee sa {\it suite caract\'eristique relative} (\`a $X$), ou s.c.r. . On a alors :
 
3. $n_i\ge m_i$. De plus, la suite caract\'eristique absolue $(m_{j}')¥$ de $Y$ est form\'ee par les entiers $m_0,\dots,n_0-1,\dots,m_{d-1},\dots,n_{d-1}-1$.
Si $c_{i}:=card\{j/m_{j}\le i\}$, $d_{i}:=card\{j/n_{j}\le i\}$, et $l_{i}':=card\{j/m_{j}'= i\}$, on a $l_{i}'=c_{i}-d_{i}$ v\'erifie la condition de croissance de
Macaulay.

4. Si $Y$, de s.c.r. $(n_{i})$, et la section $H_{s}¥$ de $X$ par une hypersurface de degr\'e $s$ sont li\'es dans une hypersurface ACM $Y'$ de $X$,
 de s.c.r. $(n_{i})'$, on a $n_{i}'=n_{i}+s$.
 
5. Si $Y$ est contenue dans $Y'$ sur $X$, alors on a : 
$n_i\le n_i'$.  En particulier, $n_{i}\ge m_{i}$, et si $X$ est irr\'eductible, $n_{i}\le n_{0}+m_{i}$.

6. Si $X$ est irr\'eductible,  on a aussi: $n_{i+1}\le n_{i}+1$¥.

7. Supposons $X$ de Gorenstein. 
Alors, pour deux hypersurfaces ACM $Y,Y'\subset X$ li\'ees dans une section $X\cap H$ de $X$ avec une hypersurface $H$ de degr\'e $s$, on a
 $n_{i}'+n_{{d-1-i}}=m_{d-1}+s$.

8. Soit $d,\alpha\in\N$, et $\alpha:=sd-r, r<d$. Soit $\Delta$ le r\'esiduel d'une intersection compl\`ete $(1,r)¥$ dans une intersection compl\`ete $(s,d)¥$.
Soit $X$ une hypersurface irr\'eductible. 
Alors, pour toute autre hypersurface ACM $Y\subset X$, on a $(\forall l\ge 0) \phi_Y(l)\ge \phi_{\Delta}(l)$. De plus 
$(\forall l\ge 0) \phi_Y(l)=\phi_{\Delta}(l)$ ssi $Y$ est comme $\Delta$ le r\'esiduel d'une intersection compl\`ete $(1,r)$
 dans une intersection compl\`ete $(s,d)$.

9. Pour une hypersurface irr\'eductible $X$¥de degr\'e $d$ g\'en\'erale, les conditions n\'ecessaires $n_i\ge i, n_{i+1}\le n_{i}+1$ ne sont pas
 suffisantes pour l'existence d'un $Y$  sur $X$ r\'ealisant $(n_i)$. Mais si $X$ est plane, elles le sont.

10. 
Soit $X$ une courbe de Gorenstein, et $Y$ localement principal sur $X$. Alors, on a :
$h^0({\bf O}_{X}(Y))=1+deg(Y)-\phi_{Y}(m_{d-1}-2)$. 
Si $X$ est une courbe plane, on r\'e\'etablit \`a partir de l\`a la description g\'eom\'etrique des syst\`emes lin\'eaires de dimension maximale sur $X$.
On discute pour terminer la question suivante. Pour tout degr\'e $s$, les hypersurfaces de degr\'e $s$ d\'eterminent sur $X$ un syst\`eme lin\'eaire complet 
de degr\'e $sd$; soit $r(sd)$ sa dimension. 
Peut-on montrer, comme pour les courbes planes, que r\'eciproquement, pour $s\le m_{d-1}-2$, toute s\'erie lin\'eaire 
de deg\'e $sd$ et de dimension $r(sd)¥$ est d\'etermin\'ee par les sections de $X$ par une hypersurface de degr\'e $s$ ?

Certaines d\'efinitions et d\'emonstrations auraient pu \^etre omises dans ce qui suit, mais nous avons pr\'ef\'er\'e supposer le minimum de connaissances, 
ce qui permet de nous adresser \`a un plus large public.
 
 {\bf Mots cl\'es : Fonction de Hilbert, intersection compl\`ete, groupes de points, syst\`emes  lin\'eaires.}

 \section{Rappels}
 
 Soit $k$ un corps alg\'ebriquement clos de caract\'eristique nulle. $k^{n+1}$ est naturellement muni d'une structure de sch\'ema, not\'e $A^{n+1}$, d'anneau $A_n:=k[X_0,\dots,X_n]$.
 Etant donn\'e un id\'eal $I$ de $A_n$ {\it homog\`ene} (i.e. engendr\'e par des polyn\^omes homog\`enes), on lui associe un {\it c\^one} $X:=C(I)$, sous-sch\'ema de $A^{n+1}$ d\'efini par:
 
 i) Son {\it support}, l'ensemble $S(I)$ des points de $k^{n+1}$ sur lesquels tous les polyn\^{o}mes de $I$ s'annulent. 
 
 ii) Son {\it anneau}, l'anneau (gradu\'e) $A_X:=A_{n}/I=\oplus_l{A_X(l)}$.
  La {\it dimension} du c\^one $X$ (ou de son anneau $A_X$) est la longueur de la plus longue suite strictement croissante d'id\'eaux premiers de $A_{X}$. La dimension est nulle ssi $S(I)=\{0\}$.
  
   Par convention, le c\^one vide $\emptyset$, associ\'e \`a l'anneau nul, a toutes les dimensions. 
   
   Pour un id\'eal homog\`ene $I$ de $A_n$, consid\'erons une d\'ecomposition primaire (homog\`ene) : $I=\cap_{i=1}^k{Q_i}$. Les radicaux 
   $P_i:=\sqrt{Q_i}$ sont bien d\'etermin\'es par $I$; ce sont ses {\it id\'eaux premiers associ\'es}.  Les {\it composantes irr\'eductibles} sont les $C(P_i)$.
    Les composantes primaires $Q_j$ dont l'id\'eal premier associ\'e $P_i=\sqrt{Q_i}$ est {\it minimal},
     i.e. ne contenant strictement aucun $P_j (j\not=i)$ sont bien d\'etermin\'ees par $I$. 
     Les autres composantes $Q_i$ ne sont pas d\'etermin\'ees univoquement par $I$; les composantes irr\'eductibles $C(P_i)$
      correspondantes sont dites {\it immerg\'ees}. 
   
   Soit $I(l):=I\cap A_n(l)$. Si $I(l)=(A_n)(l)$ pour $l$ grand, on dira que $I$ est {\it irrelevant}.
    Il revient au m\^eme de dire que sa racine est l'id\'eal maximal $(X_0,\dots,X_n)$. Sinon, $I$ est {\it relevant}. 
    $I$ est irrelevant ssi $S(I)=\{0\}$. On d\'efinit $I_{sat}$ comme l'ensemble des polyn\^omes $P$ tels que pour un certain entier $m$, 
     la multiplication de $P$ par un polyn\^ome homog\`ene quelconque de degr\'e $m$ appartient \`a $I$. 
     Alors $I_{sat}$ est aussi l'intersection des composantes relevantes d'une d\'ecomposition primaire. 
     $I=I_{sat}$ \'equivaut donc \`a l'existence d'une d\'ecomposition primaire sans composante irrelevante. 
     On dit alors que $I$ est {\it satur\'e}.
      Le c\^one $C(I)$ (resp. l'anneau $A/I$) est dit {\it projectif} si $I$ est satur\'e. 
      A un id\'eal homog\`ene $I$, on associe un faisceau $\bf I$, dont la fibre en $x\in \P_n$ est l'ensemble des fonctions rationnelles $f:=g/h, h(x)\not= 0$.
       Ce faisceau d\'efinit un sous-sch\'ema de $\P_n$, ou un {\it sous-sch\'ema projectif}.
        On le note $Proj(A_n/I)$, ou $V(I)$. Deux id\'eaux $I$ et $J$ de $A_n$ d\'efinissent donc le m\^eme sous-sch\'ema projectif
         ssi $I_{sat}=J_{sat}$, ou encore $I_l=J_l$ pour $l>>0$.  
         L'application $I\mapsto Proj(A_n/I)$ est donc bijective entre les id\'eaux homog\`enes satur\'es, 
         et les sous-sch\'emas projectifs de $\P_n$ non vides. 
         On note son inverse $X\mapsto I_X$.
          On note $A_X=A_n/I_X$ l'{\it anneau gradu\'e} (projetant) du sous-sch\'ema projectif $X$. 
          La dimension du sous-sch\'ema projectif $X$ est celle de $A_X$, moins $1$. 
          On identifiera souvent par la suite le c\^one $C(I)$ avec l'anneau gradu\'e $A_n/I$ qui le d\'efinit, 
          et avec le sous-sch\'ema projectif $V(I)$ lorque l'id\'eal est satur\'e.
           On dira que $X$ est de {\it dimension pure} $m$
            si toutes ses composantes irr\'eductibles sont de dimension $m$. 
            On appelle {\it courbe alg\'ebrique} (projective) un sous-sch\'ema projectif de dimension pure $1$, et {\it groupe de points} 
            un sous-sch\'ema projectif de dimension pure $0$. 
   
   Etant donn\'e $B$ une $k-$alg\`ebre gradu\'ee de type fini, on note $B(l)\subset B$ le $k-$sous-espace vectoriel des \'el\'ements de degr\'e $l$, 
   et son rang sur $k$,$\phi_B(l):=rg_k(B(l))$, sa {\it fonction de Hilbert}; pour $B=A_X=A_n/I_X$, on note $\phi_X$ sa fonction de Hilbert. 
   On note $B[m]$ l'alg\`ebre gradu\'ee d\'efinie par la graduation $B[m](l):=B(m+l)$.
   
   Si $X$ et $Y$ sont deux c\^ones de $A^{n+1}$, 
   d\'efinis respectivement par les id\'eaux $I_X$ et $I_Y$, on d\'efinit le c\^one $X\cap Y$ par l'id\'eal $I_X+I_Y$, et le c\^one $X\cup Y$ 
   par l'id\'eal $I_X\cap I_Y$. Si $X$ et $Y$ sont projectifs, $X\cup Y$ l'est aussi; mais $X\cap Y$ en g\'en\'eral ne l'est pas. 
   Si $X$ et $Y$ sont des sous-sch\'emas projectifs, on a donc $I_{X\cap Y}={(I_X\cap I_Y)}_{sat}$.
    Soit $X$ un c\^one de $A^{n+1}$, d'anneau $A_X$. 
    Un polyn\^ome homog\`ene $h\in A_X$ est {\it r\'egulier}
     si la multiplication par $h$ est injective dans $A_X$. 
     Pour qu'il existe un \'el\'ement r\'egulier, 
     il faut et il suffit que $X$ soit projectif. 
     Dans ce cas, si $H:=C((h))$, on a $A_{X\cap H}=A_X/hA_X$. 
     Une suite $(f_1,\dots,f_r)$ de polyn\^omes homog\`enes est {\it r\'eguli\`ere} dans $A_X$
      si pour tout $i,1\le i\le r$, la multiplication par $f_i$ est injective dans $A_X/{(f_1,\dots,f_{i-1})A_X}$.
   
    \begin{defi}
     Le c\^one $X$ (resp. son anneau $A_X$, resp. le sous-sch\'ema projectif $X':=Proj(A_X)$) est dit {\it arithm\'etiquement Cohen-Macaulay} 
     (on notera ACM par la suite) si il existe une suite r\'eguli\`ere de $m:=dim(X)=dim(X')+1$ polyn\^omes homog\`enes dans $A_X$.\end{defi}
    
    Alors,
     soit $X$ un c\^one projectif de dimension $m\ge 1$ et $h$ un \'el\'ement r\'egulier.
      Alors, $X$ est ACM ssi $X\cap C((h))$ est ACM de dimension $m-1$.
       En particulier, soit $X$ un c\^one ACM de dimension $m$. 
       Alors, toute suite r\'eguli\`ere non prolongeable est de longueur $m$. 
       Pour qu'une suite $(h_1,\dots,h_s)$ de $s\le m$ polyn\^omes homog\`enes forme une suite r\'eguli\`ere,
        il faut et il suffit que le support de $X\cap C((h_1,\dots, h_s))$ soit de dimension (pure) $m-s$.  En particulier,
         $m$ formes lin\'eaires g\'en\'eriques forment une suite r\'eguli\`ere. 
 
  \section{Suite caract\'eristique d'un sous-sch\'ema projectif ACM}
  
  Soit $X$ un c\^one ACM de dimension pure $m+1\ge 0$, $Y_0,\dots,Y_m$ une suite r\'eguli\`ere de formes lin\'eaires pour $A_X$. 
  Soit $Z_{1},\dots, Z_{p} (p:=n-m)$ un syst\`eme de coordon\'ees homog\`enes compl\'ementaire. 
  
  \begin{lem} \label{libre} Le $R_m-$module $A_X$ (o\`u $R_m:=k[Y_0,\dots,Y_m]$) admet une base $e_{0},\dots,e_{d-1}$, 
  form\'ee de mon\^omes de $k[Z_{1},\dots,Z_{p}]$ :
   $A_X\simeq\oplus_{i=0}^{d-1}{R_m[-m_i]}$, $0= m_0\le \dots\le m_{d-1}$, $m_{i}=deg(e_{i})$.
  \end{lem}
  
  \dem
   La d\'emonstration se fait par r\'ecurrence sur $m$. On commence par $m=-1$ ($R_{-1}:=k$).
   $X$ a alors pour support $\{0\}$, et $A_X(l)=0$ pour $l$ assez grand. $A_X$ est donc un $k-$espace vectoriel de dimension finie, soit $d$.
    Si $y_0,\dots,y_{l}$ est un syst\`eme de g\'en\'erateurs form\'e de mon\^omes de $k[Z_{1},\dots,Z_{p}]$, 
    on peut en extraire une base pour $A_X$, soit par exemple $y_0,\dots,y_{d}$; si on pose $m_i=deg(y_i)$, on a $A_X\simeq \oplus_{i=0}^d{k[-m_i]}$. 
  Supposons montr\'e le lemme pour $m-1$. Soit $X$ de dimension $m+1$, et $Y_0,\dots,Y_m$ une suite r\'eguli\`ere de formes lin\'eaires dans $A_X$.
   Alors la multiplication par $Y_{m}$ est injective dans $A_X$, et si $X'=X\cap\{Y_{m}=0\}$, le fait que $X$ est ACM implique que $X'$ est ACM de dimension $m$. Consid\'erons $R_{m-1}=k[Y_0,\dots,Y_{m-1}]$; d'apr\`es l'hypoth\`ese de r\'ecurrence, $A_{X'}$ admet une base $e_0,\dots,e_{d-1}$ sur $R_{m-1}$. Soit $e_0',\dots,e_{d-1}'$ dans $A_X$ tels que $e_i=e_i' mod (Y_{m})$ dans $A_{X'}$.
Alors $e_0',\dots,e_{d-1}'$ forment une base du $R_m-$module $A_X$.
  Supposons en effet donn\'ee une relation $\sum_{i=0}^n{a_i' e_i'=0}$. 
  Alors, prenons-la modulo $Y_{m}$; on obtient $\sum_{i=0}^{d-1}{a_i e_i}=0$, avec $a_i=a_i' mod Y_{m}$. 
  On en d\'eduit $a_i=0$, donc $a_i'$ est multiple de $Y_{m}$ : $a_i'=Y_{m} b_i'$. 
  Comme la multiplication par $Y_{m}$ est injective, on obtient $\sum_{i=0}^{d-1}{b_i'e_i'}=0$; 
  mais alors $b_i'$ est comme $a_i'$ multiple de $Y_{m}$. De proche en proche, on voit que $a_i'$ est multiple de $Y_{m}^j$ pour tout $j\in \N$, 
  ce qui implique $a_i'=0$ : le syst\`eme $\{e_0',\dots,e_{d-1}'\}$ est libre. Montrons maintenant qu'il est g\'en\'erateur. 
  Soit $M'$ le sous-$R_m-$module de $A_X$ engendr\'e par $e_0',\dots,e_{d-1}'$.
   Soit $x\in M'$; on a $x  mod  {Y_{m}}=\sum_{i=0}^{d-1}{a_i e_i}$ dans $A_{X'}$, 
   d'o\`u $x-\sum_{i=0}^{d-1}{a_i'e_i'}\in Y_{m}A_X$. On a donc $A_X/M'\subset Y_{m}A_X/M'$. 
   Comme $A_X/M'$ est gradu\'e, on en d\'eduit $A_X/M'=0$.
   \CQFD
  
  \begin{lem} Soit $(m_i)_{0\le i\le d-1}$ et $(m_i')_{0\le i\le d'-1}$ deux suites croissantes. Si $\sum_{i=0}^{d-1}{C_{m+l-m_i}^m}=\sum_{i=0}^{d'-1}{C_{m+l-m_i'}^m}$ pour tout entier $l$, alors $d=d'$, et $m_i=m_i'$ pour $0\le i\le d-1$.\label{uni}
  \end{lem}
  
  \dem  D\'efinissons pour une fonction $\phi : \Z\to \Z$ sa diff\'erence $\Delta\phi(l):=\phi(l)-\phi(l-1)$, 
  et posons la d\'efinition r\'ecursive $\Delta^{n+1}\phi=\Delta\Delta^n\phi$.
   Pour $\phi(l):= \sum_{i=0}^{d-1}{C_{m+l-m_i}^m}$, on v\'erifie que $\Delta^{m+1}\phi(l)$ est \'egal au nombre de $m_i$ \'egaux \`a $l$. 
   \CQFD
   
   On voit donc que la suite $(m_i)_{0\le i\le d-1}$ d\'efinie ci-dessus ne d\'epend que de la fonction de Hilbert 
   $\phi_X$ de $A_X\simeq \oplus_{i=0}^{d-1}{R_{m}¥[-m_i]}$ et la caract\'erise. Si l'on d\'efinit le degr\'e de $X$ 
   comme le coefficient de $l^m/m!$ dans son polyn\^ome de Hilbert $P_X(l)$, on voit que $d$ est \'egal au degr\'e de $X$. 
   
   \begin{defi} La suite $(m_i)_{0\le i\le d-1}$, qui caract\'erise $\phi_X$, est la {\it suite caract\'eristique} (absolue) de $X$.
   \end{defi}
   
   On pose $l_{i}:=card\{j/m_{j}=i\}$. 
    Soit $c>0$ un entier positif. On d\'efinit la $d-$i\`eme repr\'esenta\-tion de Macaulay de $c$
comme l'unique \'ecriture de $c$ de la forme $c=C_{k_d}^{d}+C_{k_{d-1}}^{d-1}\dots+C_{k_\delta}^{\delta}$, avec
$k_d>k_{d-1}>\dots>k_\delta\ge\delta>0$. On d\'efinit de plus
$c^{<d>}:=C_{k_d+1}^{d+1}+C_{k_{d-1}+1}^{d}\dots+C_{k_\delta+1}^{\delta+1}$. 
La suite $(a_l)_{l\in\N}$ s'appelle une {\it $0-$suite} si pour tout $l\in N$, $a_{l+1}\le a_l^{<l>}$.
On sait alors (cf. par
exemple \cite{Gotzmann}) que pour tout c\^one $X$, la suite $a_l:=\phi_{X}(l)$ satisfait \`a la propri\'et\'e d'\^etre une $0-$suite.
Si $X$ est de Cohen-Macaulay de dimension $m+1$, il
en est donc de m\^eme des suites $a_l^s:=\Delta^s\phi_X(l)$
 pour tout
$s,0\le s\le m+1$. 

Comme $l_{i}=\Delta^{m+1}¥\phi_X(l)$, on en d\'eduit:
\begin{lem}
La suite $l_{i}$ est une $0-$suite.
\end{lem}

On voit en particulier que de $l_{1}=p$, on tire $l_{j}\le C_{p-1+j}^{j}$. On voit de la d\'emonstration du premier lemme que le plus petit degr\'e d'une hypersurface 
contenant $X$ est \'egal \`a $\min\{j, l_{j}< C_{p-1+j}^{j}\}$.
   
   Soit $H:=\{h=0\}$ un hyperplan coupant $X$ proprement. La d\'emonstration du lemme \ref{libre} nous montre que la suite de $X\cap H$ est la m\^eme que celle de $X$. On a aussi $\phi_{X\cap H}=\Delta\phi_X$. 
   
\begin{lem} Soit $Y_0,\dots,Y_m$ une suite r\'eguli\`ere pour $A_X$. Soit $Y$ une combinaison lin\'eaire des coordonn\'ees homog\`enes $X_i$, $Y:=\sum_{i=0}^n{a_i X_i}$, qui n'est pas combinaison $k-$lin\'eaire des $Y_j$. Alors soit $\pi_Y : \P_n\to \P_{m+1},(X_0:\dots:X_n)\mapsto(Y_0:\dots:Y_m:Y)$. On peut choisir $Y$ de sorte que $X'=\pi_{Y}(X)$ soit une hypersurface de degr\'e $d$ ne passant pas par $(0:\dots:0:1)$; alors $1,\dots,Y^{d-1}$ forment 
une base du sous-$R_m-$module $A_{X'}\subset A_X$. De plus, si $m:=dim(X)\ge 1$, il en est ainsi pour tout choix de $Y$.\end{lem}

\dem Par hypoth\`ese, $\{Y_0=0,\dots, Y_m=0\}$ ne rencontre par $X$. 
Par cons\'equent, le point $(0:\dots:0:1)$ n'appartient pas au support de $X'$. Si $m\ge 1$, le degr\'e n'est pas chang\'e. On en d\'eduit que $1,\dots,Y^{d-1}$ est
 une base du $R_m-$module $A_{X'}$. Si $m=0$ ($X$ est un groupe de points), on peut choisir $Y$
  de fa\c{c}on \`a \'eviter les points align\'es, de sorte \`a avoir $deg(X')=deg(X)$.
\CQFD

Les suites $(m_{i})$ et $(l_{i}:=card\{j/m_{j}=i\}$ v\'erifient les propri\'et\'es suivantes:

\begin{thm} $l_0=1$, $l_1=p$, o\`u $p$ est la codimension de $X$, suppos\'e non-d\'eg\'en\'er\'e (i.e. non contenu dans un hyperplan).
 De plus, l'ensemble des indices $j$ tels que $l_j\not=0$ est {\it connexe}. 
Enfin, supposons $m\ge 1$. Alors, si $l_i=1$, et $l_{i+1}\not= 0$, on a $l_{i+1}\ge p$.
\end{thm}  
  \dem
Comme $\phi_X(0)=1$ et $\phi_{X}(-1)=0$, il y a exactement un $m_i$ \'egal\`a $0$, i.e. $m_0=0$ et $m_1>0$.
Soit $M_1,\dots,M_l$ les \'el\'ements de degr\'e $1$.
Consid\'erons une base compl\'ementaire $Y_{m+1},\dots, Y_n$ de $A_n(1)$. Ecrivons $Y_{m+j}=Z_j+a_1^j M_1+\dots+a_l^j M_l$, pour $j=1,\dots,n-m=p$, et $a_i^j\in k, Z_j\in R_m(1)$. Si on avait $l< p$, on en d\'eduirait une relation lin\'eaire entre les $Y_i$, ce qui est impossible puisque $X$ est non d\'eg\'en\'er\'e. 
D'autre part, $1,M_1,\dots,M_l$ sont des $R_m-$combinaisons ind\'ependantes de $1,Y_{m+1},\dots,Y_{m+p}$. Ce serait impossible si $l>p$. On a donc $l=p$.
Supposons $l_i=1$. Soit $M_j$ l'\'el\'ement de degr\'e $i$. Alors supposons $l_{i+1}<p$. 
Alors on pourrait trouver une
combinaison
$Y=a_1Y_{m+1}+\dots+a_{s}Y_{n}$ $k-$lin\'eaire des $Y_{m+1},\dots,Y_n$ telle que $Y M_j$ appartienne au
$R_{m}¥-$module engendr\'e par $M_0, \dots, M_j$. Mais alors pour tout $l<i$, $Y M_l$ appartient aussi
\`a ce $R_m-$module, puisque $M_j$ est le seul g\'en\'erateur en degr\'e $i$. Le th\'eor\`eme de
Cayley-Hamilton nous donne alors un polyn\^ome de degr\'e $j+1$ en $Y$ \`a coefficients dans $R_m$ qui est nul.
Mais le lemme pr\'ec\'edent nous dit que $1,Y,\dots,Y^{d-1}$ sont $R_m-$lin\'eairement ind\'ependants, si $m\ge 1$.
Si $m\ge 1$, on en d\'eduit $j=d-1$, i.e. si $l_{i+1}\not=0$, $l_{i+1}\ge p$.

Il y a deux mani\`eres possibles pour montrer $m_{i+1}\le m_{i}+1$. La premi\`ere est de dire que $l_{i}$ est une $0-$suite. Pour la deuxi\`eme, consid\'erons la multiplication par une forme lin\'eaire $Y$ dans $A_X$. Soit $M_0,\dots,M_{d-1}$ une base du 
$R_m-$module $A_X$, avec $m_i=deg(M_i)$. On a :
$Y M_i=\sum_{j=0}^{d-1}{m_{ij}M_j}$, avec $deg(m_{ij})=m_i-m_j+1$. 
Supposons $m_{i+1}>m_i+1$, pour un $i,0\le i\le d-2$. Alors on aurait $m_{kj}=0$ pour $k\le i,j\ge
i+1$, pour des raisons de degr\'es. On en d\'eduit que $Y$ d\'etermine un endomorphisme du $R_m-$module
engendr\'e par $M_0,\dots, M_i$, de matrice $M'$, sous-matrice de $M=(m_{ij})$. D'apr\`es
Cayley-Hamilton, on en d\'eduit $\det(M'-Y Id_{i+1})M_0=0$ dans $A_X$. Mais $M_0$ est de
degr\'e $0$, donc on peut supposer $M_0=1$. On en d\'eduit $P(Y_0,\dots,Y_m,Y)=\det(M'-Y Id_{i+1})=0$ dans
$A_X$. Mais $P$ est un polyn\^ome de degr\'e $i+1<d$. Cela est impossible si l'on choisit $Y$ d'apr\`es le lemme
pr\'ec\'edent tel que $deg(X)=deg(X')$. On a donc $m_{i+1}\le m_i+1$. 
\CQFD

\begin{cor} Supposons $X$ de dimension $m\ge 1$ et de codimension $p\ge 2$. Alors : $m_{d-1}\le [(2d-1)/3]$.
\end{cor}

\dem On a en effet $m_{i+3}-m_i\le 2$. Sinon, on aurait $m_{i+3}=m_i+3$, donc $m_{i+1}=m_i+1$ et $m_{i+2}=m_i+2$,
ce qui est impossible d'apr\`es la proposition pr\'ec\'edente.
\CQFD

\begin{cor} Soit $X$ un groupe de points.
Soit $r$ le nombre maximum de points align\'es de $X$. Alors $m_{d-1}\le r+[(2d-2r-1)/3]$.
\end{cor} 

\dem On choisit $Y_0$ tel que $\{Y_0=0\}$ ne rencontre aucun des points d'intersection de deux droites distinctes
joignant deux points de $X$. Alors la projection \`a partir d'un point de $\{Y_0=0\}$ applique $X$ sur un groupe
de points de degr\'e
$\ge d-r+1$. L'in\'egalit\'e $m_{i+3}-m_i\le 2$ reste donc valable tant que $i\le d-r$, d'apr\`es la d\'emonstration pr\'ec\'edente. Pour $i\ge r$, on a
$m_{i+1}\le m_i+1$.
\CQFD

{{\bf Remarque.} Lorsque que $X$ est un groupe de points, 
on d\'efinit classiquement son {\it indice de s\'eparation} $e(X)$ comme l'entier $\max\{l,h^1({\bf I}_X(l)\not= 0\}$
, entier au-del\`a duquel $\phi_X$ est constante. On voit alors facilement $e(X)=m_{d-1}-2$.}

\begin{cor} Si $X\subset \P_n$ est un groupe de points, $\phi_X$ est strictement croissante, jusqu'\`a \^etre
constante.
\end{cor}

\dem On a en effet $\phi_X(l)=\sum_{i=0}^l{l_i}$. 
Si $\phi_X(l+1)=\phi_X(l)$,a lors $l_{l+1}=0$, et $\phi_X(j)=\phi_X(l)$ pour tout $j\ge l$. 
\CQFD

Pour $X$ est une intersection compl\`ete g\'en\'erique, de dimension $m$, on peut calculer explicitement les mon\^omes $(e_{i})$ g\'en\'erateurs
 du $R_m$-module $A_X$. 
Soit $L:=\{Y_0 =\dots = Y_m\}\subset \P_{n}$ une sous-vari\'et\'e
lin\'eaire de codimension $m+1$. Soit $Z_{1},\dots, Z_{p}¥$ des coordonn\'ees homog\`enes
compl\'ementaires. Alors:

\begin{lem} Pour une intersection compl\`ete g\'en\'erique $X$ (ne rencontrant pas $L$), les mon\^omes
$Z_{1}^{i_1}\dots Z_{p}^{i_p},0\le i_1\le d_1-1,\dots,0\le i_p\le d_p-1$ forment 
une base du $R_{m}-$module $A_X$. 
En particulier, pour un tel  $X$ g\'en\'erique, les \'equations de $X$ peuvent s'\'ecrire de mani\`ere unique sous la forme:
$$Z_{1}^{d_1}= \sum_{I=(i_0,\dots,i_n), i_1+\dots+i_n=d_1, i_{m+1}<d_1,\dots, i_n<d_p }{c_{I}^1 Y^I},$$
$$\dots$$
$$Y_{n}^{d_p}=\sum_{I=(i_0,\dots,i_n),i_1+\dots+i_n=d_p,i_{m+1}<d_1,\dots,i_n<d_p}{c_{I}^p Y^I},$$
o\`u $Y^I:=Y_{0}^{i_{1}}\dots Y_{m}^{i_{m}} Z_{1}^{i_{m+1}}\dots Z_{p}^{i_{n}}$ et $c_{I}^s\in R_{m}$.
\end{lem}

\dem 
On se ram\`ene au cas $m=-1$, i.e.: 
pour des coordonn\'ees homog\`enes $Y_j:=\sum_{i=0}^{p-1}{a_{ij}X_i} (0\le j\le p-1)$ g\'en\'eriques, 
les mon\^omes
$$Y^I:=Y_0^{i_0}\dots Y_{p-1}^{i_{p-1}}, i_0<d_1,\dots,i_{p-1}<d_{p}$$ forment
 une base du
$k-$espace vectoriel
$k[X_0,\dots,X_{p-1}]/(P_1,\dots,P_p)$, lorsque les polyn\^omes homog\`enes $P_i$, de degr\'es respectifs $d_i$,
n'ont pour z\'ero commun que l'origine de $k^{p}$.
\CQFD

On retrouve en particulier les entiers $m_{I}=i_{1}+\dots+i_{p}$, pour les intersections compl\`etes.

\begin{lem} Si $X'\subset X$ est un sous-sch\'ema de $X$, alors pour tout $m\ge 0$,
$\Delta^j(\phi_X-\phi_{X'})\ge 0$ pour tout
$j\le h(X)$ (o\`u $h(X)$ est la longueur maximale d'une suite r\'eguli\`ere de $A_X$).
\end{lem}

\dem  On fait une r\'ecurrence sur $h(X)$. Pour $h(X)=0$, l'in\'egalit\'e $\phi_X\ge \phi_{X'}$ provient de la surjectivit\'e de $A_{X}\to A_{X'}$. Supposons montr\'e le lemme pour $h(X)<s$; supposons $h(X)=s$. Alors soit $h$ un \'el\'ement r\'egulier de $A_X$. Alors l'hyperplan $H=\{h=0\}$ coupant $X$ proprement, et $\Delta\phi_X=\phi_{X\cap H}$. Comme $X'\cap H\subset
X\cap H$, on en d\'eduit $\Delta\phi_X\ge \phi_{X'\cap H}\ge \Delta\phi_{X'}$, donc
$\Delta(\phi_X-\phi_{X'})\ge 0$. De plus, $h(X\cap H)<s$. L'hypoth\`ese de r\'ecurrence nous permet de conclure.
\CQFD

\begin{cor} Si $X$ et $X'\subset X$ sont deux sous-sch\'emas ACM projectifs de dimension $m$, la suite carac\-t\'eristique de $X'$ est
"incluse" dans celle de $X$.
\end{cor}

\dem  En appliquant le lemme ci-dessus, avec $h(X)=m+1$, on voit d'apr\`es ce qui pr\'ec\`ede que le nombre $l_j'$ de
$m_i'$ \'egaux \`a $j$ est inf\'erieur au nombre $l_j$ de $m_i$ \'egaux \`a $j$, pour tout $j$.
\CQFD

\begin{lem}\label{Gore} Supposons $X$ de Gorenstein. Alors on a : $m_i+m_{d-1-i}=m_{d-1}$.
\end{lem}
\dem
Soit $\omega_{A_X}$ le module dualisant de $A_{X}$. On a, comme $A_{X}¥$ est un $R_m-$module libre, $\omega_{A_X}\simeq Hom_{R_m}(A_X,\omega_{R_m})$. Comme $\omega_{R_m}\simeq
R_m[-m-1]$, on en d\'eduit $\omega_{{A_{X}}}\simeq \oplus_{{i=0}}^{d-1}{R_m[m_{i}-m-1]}$. D'autre part, comme $X$ est de Gorenstein, on a $\omega_{A_X}\simeq
A_X[s]$ pour un certain entier $s$. 

On en d\'eduit que les deux suites croissantes $(m_i-m-1)_i$ et
$(s-m_{d-1-i})_i$ co\"\i ncident, et donc : $s=m_{d-1}-m-1$, et $m_i+m_{d-1-i}=m_{d-1}$.
\CQFD

\section{Suite caract\'eristique relative}

Soit maintenant $Y\subset X$ une hypersurface ACM de $X$. On va introduire la {\it suite
caract\'eristique relative} (\`a
$X$) de $Y$, pour profiter de l'information que $Y$ est contenu dans $X$.

Soit $Y\subset X$, avec $X$ ACM de dimension pure $m\ge 1$ et $Y$ ACM de dimension pure $m-1$ dans $X$. Le morphisme surjectif $A_X\to A_Y$ a un noyau $I_{Y/X}\simeq I_Y/I_X$. Soit $Y_0,\dots,Y_m$ $m$ formes lin\'eaires d\'efinissant un sous-espace projectif ne rencontrant pas $X$. On pose encore $R_m:=k[Y_{0},\dots,Y_{m}]$. 

\begin{lem} Pour $Y\subset X$, $I_{Y/X}$ est un $R_m-$module libre gradu\'e.
\end{lem}
\dem Consid\'erons la suite exacte $0\to I_{Y/X}\to A_X\to A_Y\to 0$.
 On a vu que $A_X$ est un $R_m-$module libre gradu\'e.
  Par ailleurs, $Y_0=\dots=Y_m=0$, ne rencontrant pas $X$, ne rencontre pas non plus $Y$. 
  Le th\'eor\`eme des syzygies gradu\'e (cf. Appendice) permet de conclure que $I_{{Y/X}}$ est aussi un $R_{m}-$module libre.
  \CQFD

\begin{defi} La suite $(n_i)$ est la suite caract\'eristique relative (\`a $X$) de $Y$. On notera par la suite s.c.r. pour suite caract\'eristique relative.\end{defi}

Lorsque $Y$ est de codimension deux et $X$ est une hypersurface de degr\'e minimal qui le contient,
 on retrouve le {\it caract\`ere num\'erique} de $Y$, introduit par Gruson et Peskine dans (\cite{Peskine}).

\begin{lem} Soit $X$ et $Y\subset X$ deux sous-sch\'emas projectifs.
 Soit $h$ une forme lin\'eaire qui d\'etermine une multiplication injective dans $A_X$ et dans $A_{Y}$, et $H:=\{h=0\}$.
On a $I_{{Y\cap H}/{X\cap H}}\simeq I_{Y/X}/hI_{Y/X}$.\end{lem}

\dem Appliquons \`a la suite exacte $0\to I_{Y/X}\to A_X\to A_Y\to 0$ le foncteur $\otimes_{A_n}A_n/hA_n$. On obtient, comme $M\otimes_A A/I\simeq M/IM$ et que le foncteur est exact \`a droite, la suite exacte : $I_{Y/X}/hI_{Y/X}\to A_X/hA_X\simeq A_{X\cap H}\to A_Y/hA_Y\simeq A_{Y\cap H}\to 0$. On veut montrer que la premi\`ere fl\`eche de cette suite est injective. Pour cela, remarquons que la suite exacte est en particulier une suite exacte de $k-$espaces vectoriels. Puisque la multiplication par $h$ est injective dans $A_X$, elle l'est aussi dans $I_{Y/X}$, et donc $\phi_{I_{Y/X}/hI_{Y/X}}=\Delta \phi_{I_{Y/X}}$. Mais $\phi_{I_{Y/X}}=\phi_{A_X}-\phi_{A_Y}$. On en d\'eduit
$\phi(I_{Y/X}/hI_{Y/X})=\Delta\phi_{X}-\Delta\phi_{Y}=\phi_{X\cap H}-\phi_{Y\cap H}$, ce qui nous
montre que dans la suite exacte de $k-$espaces vectoriels 
$I_{Y/X}/hI_{Y/X}\to A_{X\cap H}\to A_{Y\cap
H}\to 0$, la premi\`ere fl\`eche est injective (sinon son noyau nous conduirait \`a
une contradiction). On a donc bien une suite exacte de $A_n-$modules $0\to I_{Y/X}/hI_{Y/X}\to
A_{X\cap H}\to A_{Y\cap H}\to 0$, donc $I_{Y/X}/hI_{Y/X}\simeq I_{Y\cap H/X\cap H}$.
\CQFD

Dans ce qui suit, $X$ est un sous-sch\'ema projectif ACM, de suite caract\'eristique $(m_i)$, et $Y$ une hypersurface de $X$, de s.c.r. $(n_i)$.
On d\'eduit du lemme pr\'ec\'edent:
\begin{cor} Si l'hyperplan $H$ coupe $X$ et $Y$ proprement, la s.c.r. de $Y\cap H$ dans $X\cap H$ est la m\^eme que celle de $Y$ dans $X$.
\end{cor}

\begin{thm}
Si $Y\subset X$ ($m\ge 1$) est une hypersurface ACM de $X$, on a 
$(\forall i) n_i\ge m_i$.
De plus, si $Y'$, de s.c.r. $(n_i')$, contient $Y$, on a :
$(\forall i) n_i'\ge n_i$.
\end{thm}

\dem
Montrons $n_i'\ge n_i$.
Soit $H:=\{h=0\}$ un hyperplan coupant $X$, $Y$, et $Y'$ proprement. On a encore $Y\cap H\subset Y'\cap H$ dans
$X\cap H$.
Par ailleurs, la s.c.r. de $Y\cap H$ (resp. $Y'\cap H$) dans $X\cap H$ reste
inchang\'ee. Il suffit donc de montrer l'\'enonc\'e pour $m=0$, i.e. pour $Y$ et $Y'$ de support $\{0\}$ dans $A^{n+1}$. 

Posons $1_k(l)=(l+1-k)_+-(l-k)_+= 1$ si $l\ge k$, $0$ sinon.
$\phi_Y(l)\le \phi_{Y'}(l)$ se lit alors :

 $$\sum_{j=0}^{d-1}{1_{n_j'}(l)}\le
\sum_{i=0}^{d-1}{1_{n_j}(l)}.$$

L'in\'egalit\'e ci-dessus nous donne pour $l=n_0'$ que $n_0\le n_0'$.
Supposons qu'on ait montr\'e que pour tout $j< i,n_j\le n_j'$.

Alors pour un $j< i$, $1_{n_j'}(n_i')=1$, et $1_{n_j}(n_i')=1$ car $n_j\le n_j'\le n_i'$.

Si $n_i'<n_i$, on aurait pour $j\ge i$ $n_i'<n_i\le n_j$ donc $1_{n_j}(n_i')=0$. L'in\'egalit\'e
ci-dessus ne serait donc pas v\'erifi\'ee puisque $1_{n_i'}(n_i')=1$. On a donc $n_i'\le n_i$, pour
tout $i$.

Montrons maintenant $n_i\ge m_i$. On a une injection gradu\'ee $\oplus_i{R[-n_i]}\to \oplus_i{R[-m_i]}$. De  $\sum_{j=0}^{d-1}{1_{n_j}(l)}\le
\sum_{i=0}^{d-1}{1_{m_j}(l)},$ on d\'eduit $n_{i}\ge m_{i}$.
\CQFD

Consid\'erons deux r\'esolutions libres $0\to \oplus_i{R_m[-n_i]}\to \oplus_i{R_m[-m_i]}\to A_Y\to 0$ et
$0\to \oplus_i{R_m[-n_i']}\to \oplus_i{R_m[-m_i']}\to A_Y\to 0$ (\'eventuellement infinies) du sous-sch\'ema ACM $Y$. 
Alors, le lemme \ref{uni} nous montre que $(n_i)\oplus (m_i')=(n_i')\oplus (m_i)$, o\`u l'on d\'efinit la somme $\oplus$ des suites en additionnant les largeurs $l_{i}¥$ ¥associ\'ees. 
En particulier, si la suite $(m_i')$ (resp. $(n_i')$ est obtenue \`a partir de $(m_i)$ en supprimant certaines valeurs, la suite $(n_i')$ (resp. $(m_i')$) est obtenue \`a partir de $n_i$ (resp. de $(m_i)$) en supprimant les m\^emes valeurs. 
En particulier, on en d\'eduit la relation entre les suites caract\'eristiques absolues et relatives:

\begin{lem} 
La suite caract\'eristique absolue $(m_{i}'¥)¥$ de $Y$ est \'egale \`a : $$(m_0,m_0+1,\dots,n_0-1)\oplus (m_1,\dots,n_1-1)\oplus (m_{d-1},\dots,n_{d-1}-1).$$
\end{lem}

\dem On consid\`ere la suite exacte $0\to \oplus_{i=0}^{d-1}{R_m[-n_i]}\to \oplus_{i=0}^{d-1}{R_m[-m_i]}\to A_Y\to 0$.
On \'ecrit
$R_m=\oplus_{j=0}^\infty{R_{m-1}[-j]}$, o\`u $R_{m-1}=k[Y_0,\dots,Y_{m-1}]$. On en d\'eduit la suite exacte de $R_{m-1}-$modules :
$$0\to \oplus_{0\le i\le d-1,0\le j<\infty}{R_{m-1}[-n_i-j]}
\to \oplus_{0\le i\le d-1,0\le j<\infty}{R_{m-1}[-m_i-j]}\to
A_Y\to 0$$ 
Par ailleurs, on suppose que $\{Y_0=0,\dots,Y_{m-1}=0\}$ coupe $Y$ proprement. Alors on sait que $A_Y\simeq \oplus_{j}{R_{m-1}[-m_j']}$
 pour certains entiers $m_{j}'$¥. La remarque pr\'ec\'edant le th\'eor\`eme montre que la suite $(m_i')$ est obtenue \`a partir de la suite exacte pr\'ec\'edente
 en supprimant les doubles.
\CQFD

On peut reconstruire en sens inverse la suite $(n_{i})$ avec des suites caract\'eristiques absolues $(m_{i})$ et $(m_{i}')$ de $X$ et $Y$. 
En particulier, on voit que si $c_{i}:=card{\{j/m_{j}\le i\}}$, $d_{i}:=card{\{j/n_{j}\le i\}}$, et $l_{i}':=card\{j/m_{j}'=i\}$, 
on a $l_{i}'=c_{i}-d_{i}$. On a donc la propri\'et\'e suppl\'ementaire sur les $n_{i}$:
\begin{lem}
$c_{i}-d{i}$ est une $0-$suite.
\end{lem}

Supposons que le sous-sch\'ema $X$ est de Gorenstein, et que $Y$ et $Y'$ sont li\'ees dans la section de $X$ par une hypersurface $H$ de degr\'e $s$  (cf. Appendice).
 Alors, on peut dans ce cas relier entre elles les s.c.r. respectives $(n_i)$ et $(n_i')$ de $Y$ et de $Y'$ par la formule suivante:

\begin{thm}
$n_i+n_{d-1-i}'=m_{d-1}+s$.
\end{thm}

\dem
Soit $Y$ et $Y'$ li\'ees dans $X\cap H$, avec $H$ d\'efini par un polyn\^ome homog\`ene de degr\'e $s$.
Consid\'erons la suite exacte $0\to I_{Y/X}\to A_X\to A_{Y}\to 0$.
On lui applique le foncteur $Hom_{A_X}(\bullet,A_X)$.
D'abord $Hom_{A_X}(A_X,A_X)\simeq A_X$; l'isomorphisme associe \`a $\phi\in Hom_{A_X}(A_X,A_X)$ sa valeur
$\phi(1)$. De plus, $Hom_{A_X}(A_Y,A_X)=0$. En effet, soit $\phi\in Hom_{A_X}(A_Y,A_X)$. Comme
$A_Y=A_X/(I_{Y/X})$, on peut associer \`a $\phi$ canoniquement un morphisme $\tilde \phi\in Hom_{A_X}(A_X,A_X)$,
donc un \'el\'ement
$a=\tilde\phi(1)\in A_X$. On doit avoir $a x=0$ pour tout $x\in I_{Y/X}$. Choisissons $x$ dans $I_Y$, et $x$
n'appartenant \`a aucun des id\'eaux premiers $P_i$ associ\'es \`a $I_X$ (cela est possible d'apr\`es le lemme
d'\'evitement puisque $I_Y$ ne peut \^etre contenu dans aucun des $P_i$). Alors $a cl(x)=0$ implique $a=0$.

Consid\'erons maintenant $Hom_{A_X}(I_{Y/X},A_X)$.
La multiplication par $h$ d\'efinit un isomorphisme  de $A_{X}-$modules $I_{X\cap H/X}\simeq A_X[-s]$.
A un morphisme $\phi\in Hom_{A_X}(I_{Y/X},A_X)$ correspond donc bijectivement un morphisme de $I_{Y/X}$ dans $I_{H\cap X/X}[s]$.
Comme $Y$ et $Y'$ sont li\'es dans $X\cap H$, on a $Hom_{A_X}(I_{Y/X},I_{X\cap H/H})\simeq I_{Y'/X}$.
On a donc bien $Hom_{A_X}(I_{Y/X},A_X)\simeq I_{Y'/X}[s]$.

Enfin, on a \'evidemment $Ext^1_{A_X}(A_X,A_X)=0$, 
d'o\`u une suite exacte de $A_X-$modules : $0\to A_X\to I_{Y'/X}[s]\to Ext^1_{A_X}(A_Y,A_X)\to 0$.
On en d\'eduit une suite exacte de $R_{m}$-modules: $$0\to \oplus{R_{m}¥[-m_i]}\to \oplus{R_{m}[-n_i'+s]}\to
Ext^1_{A_X}(A_Y,A_X)\to 0$$ 

De m\^eme, appliquons \`a la suite exacte $0\to I_{Y/X}\to A_X\to A_{Y}\to 0$
le foncteur
$Hom_{R_{m}}(\bullet,R_{m})$. 

On a $Hom_{¥R_{m}}(A_Y,R_{m}¥)=0$. En effet, soit $\phi\in Hom_{R_{m}}(A_Y,R_{m}¥)$.
 Comme $A_Y=A_X/(I_{Y/X})$, on peut lui associer
un morphisme $\tilde\phi\in Hom_{R_{m}}(A_X,R_{m})$, qui est d\'etermin\'ee par ses valeurs $e_0,\dots,e_{d-1}\in R_{m}¥$ sur une
base $f_0,\dots,f_{d-1}$ du $R_{m}-$module
$A_X$. $\tilde\phi$ doit s'annuler sur $I_{Y/X}$; soit
$\alpha_i=\sum_j{a_{i,j}f_j}$ une base de $I_{Y/X}$. On a pour tout $i$ :
$\sum_{j}{\alpha_{i,j}e_j}=0$. Mais comme les $\alpha_i$ sont
$R_{m}-$ind\'ependants, il n'y a pas de relation sur $R_{m}¥$ non triviale entre les colonnes de la matrice $a_{i,j}$, donc
pour tout
$j$,
$e_j=0$.

Enfin,
$Ext^1_{R_{m}}¥(A_X,R_{m})=0$ (puisque
$A_X$ est un
$R_{m}¥-$module libre). On obtient donc la suite exacte
$$0\to
\oplus_i{R_{m}[m_i]}\to
\oplus_i{R_{m}[n_i]}\to Ext^1_{R_{m}}¥(A_Y,R_{m})\to 0$$

Mais $Ext^1_{R_{m}}¥(A_Y,\omega_{R_{m}}¥)\simeq Ext^1_{A_X}(A_Y,\omega_{A_X})\simeq \omega_{A_Y}$
comme
$R_{m}¥-$modules, avec $\omega_{R_{m}}¥=R[-m-1]$ et $\omega_{A_X}\simeq A_X[m_{d-1}-m-1]$ comme on l'a vu.

D'o\`u les suites exactes :
 $$0\to \oplus{R_{m}¥[m_{d-1}-m-1-m_i]}\to \oplus{R_{m}¥[-n_i'+s+m_{d-1}-m-1]}\to
\omega_{A_Y}\to 0$$

et : $$0\to \oplus{R_{m}¥[m_i-m-1]}\to
\oplus{R_{m}¥[n_i-m-1]}\to \omega_{A_Y}\to 0.$$

 Comme les suites $(n_i)$ et $(n_i')$ sont par d\'efinition croissantes, et comme on a vu
que les suites
$(m_{d-1}-m-1-m_{d-1-i})_i$ et
$(m_i-m-1)_i$ co\"\i ncident, on en d\'eduit l'\'egalit\'e des suites
$(m_{d-1}-m-1+s-n_{d-1-i}')_i$ et $(n_i-m-1)_i$; donc $n_i+n_{d-1-i}'=s+m_{d-1}$.
\CQFD

\begin{cor}
Soit $X$ est de Gorenstein, et $s$ le 
plus petit degr\'e d'une hypersurface contenant $Y$ sans contenir aucune composante de $X$. Alors,
$n_{i}\le s+m_{i}$. En particulier, si $X$ est irr\'eductible, on retrouve $n_i\le n_0+m_i$.
\end{cor}
\dem
Soit $H$ une hypersurface de degr\'e $s$ contenant $Y$, et qui coupe $X$ proprement.
Soit alors $Y'$ le r\'esiduel de $Y$ dans $X\cap H$. Soit $(n_i')$ sa s.c.r. .
On a d'apr\`es le th\'eor\`eme pr\'ec\'edent $n_i=s+m_{d-1}-n_{d-1-i}'$; comme $n_{d-1-i}'\ge
m_{d-1-i}$, on a : $n_i\ge s+m_{d-1}-m_{d-1-i}=s+m_i$.
Si $X$ est irr\'eductible, le premier g\'en\'erateur $\alpha_0$ de $I_{Y/X}$ d\'efinit une
hypersurface de degr\'e $n_0$ contenant $Y$ et coupant $X$ proprement, donc $n_{i}\le n_{0}+m_{i}¥$.
\CQFD

Soit $H$ une hypersurface de degr\'e $s$. Si $H_{s}:=X\cap H$ et l'hypersurface ACM $Y$¥sont li\'es dans une autre hypersurface 
ACM $Y'=Y\cup H_s$, on peut calculer la
s.c.r. $(n_{i}'¥)¥$¥ de $Y'$
 \`a partir de celle $(n_{i})$ de $Y$, par la formule suivante:

\begin{pro}
$n_i'=n_i+s$. \label{sexion}
\end{pro}

\dem
En effet, si $h$¥ est un polyn\^ome homog\`ene de degr\'e
$s$ d\'efinissant
$H_s$, la
multiplication par $h$ d\'etermine un isomorphisme de $A_{X}-$modules : $$I_{Y/X}[-s]\simeq I_{Y'/X}.$$ 
\CQFD

Le th\'eor\`eme suivant g\'en\'eralise un r\'esultat de \cite{EP}:
\begin{thm} Soit $X\subset \P_n$ une hypersurface r\'eduite, et $Y\subset X$ une hypersurface de $X$, de s.c.r.
 $(n_{i})$. Supposons
$n_{i}>n_{i-1}+1$. Alors,
$X=X'\cup X''$, avec $X'$ et $X'$ des hypersurfaces de degr\'es respectifs
$deg(X')=i$ et $deg(X'')=s:=d-i$. De plus, si l'on pose $Y'=Y\cap X'$ et $Y''=Y\cap X''$, la s.c.r. de $Y''$ dans
$X''$ est $(n_i,\dots,n_{d-1})$. Celle de $Y'$ dans $X'$ est
$(n_0-s,\dots,n_{i-1}-s)$.
\end{thm}

\dem Soit $(\alpha_i)$ une base du $R_{m}¥-$module $I_{Y/X}$. Consid\'erons la multiplication par une forme lin\'eaire
$T$ :
$T\alpha_i=\sum_{j=0}^{d-1}{t_{ij}\alpha_j}$, avec $t_{ij}\in R_{m}$. Supposons $n_{i}>n_{i-1}+1$. Alors la matrice
$(t_{kj})$ doit avoir les termes $k\le i-1,j\ge i$ nuls. En particulier, le polyn\^ome caract\'eristique de la
matrice
$(t_{ij})$ s'\'ecrit comme un produit $P(T)Q(T)$, o\`u $P(T)$ est le polyn\^ome caract\'eristique de la sous-matrice
$(t_{kj})_{0\le k,j\le i-1}$. Le sous-$R_{m}-$module engendr\'e par $\alpha_0,\dots,\alpha_{i-1}$ est stable par
multiplication par $T$. Le th\'eor\`eme de Cayley-Hamilton nous donne 
$$P(T)(\alpha_0,\dots,\alpha_{i-1})=0, P(T)Q(T)(\alpha_0,\dots,\alpha_{d-1})=0.$$

 Si on choisit un \'el\'ement $\alpha$ de $I_{Y/X}$ qui ne s'annule sur
aucune composante de $X$ (par le lemme d'\'evitement), la relation $P(T)Q(T)\alpha=0$ nous donne $P(T)Q(T)=0$.
On obtient ainsi une hypersurface de degr\'e $d$ contenant $X$; comme $X$ est r\'eduite de degr\'e $d$, c'est $X$
elle-m\^eme. Si on pose $X':=\{P=0\}$, $X'':=\{Q=0\}$, on a : $X=X'\cup X''$. 
D'autre part, comme $\alpha_0,\dots,\alpha_{d-1}$ g\'en\`erent $I_{Y/X}$, ils g\'en\`erent aussi $I_{Y''/X''}$.
Mais $\alpha_0,\dots,\alpha_{i-1}=0$ dans $A_{X''}$. Donc $\alpha_i,\dots,\alpha_{d-1}$ sont g\'en\'erateurs de
$I_{Y''/X''}$. On a donc
$I_{Y''/X''}\simeq \oplus_{j=i}^{d-1}{R[-n_j]}$.

Enfin, on a une suite exacte de $A_n-$modules, donc de 
$R_{m}-$modules : $$0\to I_{Y'/X'}[s]\to I_{Y/X}\to
I_{Y''/X''}\to 0.$$ 
En effet, soit $\alpha\in I_{Y/X}$.  Si $\alpha\in  I_{X''/X}$, alors $\alpha$ est un multiple de $Q$, et est donc
dans l'image de $I_{Y'/X'}[s]\to I_{Y/X}$. 
Cette derni\`ere application est injective, car si $Qx\in I_X$, alors
$Qx\in I_{X'}$, donc, comme $P$ et $Q$ n'ont pas de facteur commun, $x\in I_{X'}$.
 Enfin, il est \'evident que $I_{Y/X}\to
I_{Y''/X''}$ est surjective, car si $x=y mod I_{X''}$, alors $x=a_{i+1}\alpha_{i+1}+\dots+a_{d-1}\alpha_{d-1} mod
I_{X''}$, et $a_{i+1}\alpha_{i+1}+\dots+a_{d-1}\alpha_{d-1}\in I_{Y/X}$ a pour image $x$.

 On en d\'eduit $I_{Y'/X'}\simeq \oplus_{j=0}^{i-1}{R_{m}¥[n_j-s]}$. 
\CQFD

{{\bf Remarque.} La d\'emonstration pr\'ec\'edente donne explicitement les g\'en\'e\-rateurs de $I_{Y'/X'}$ et de
$I_{Y''/X''}$. De $P\alpha_j=0$ pour $j\le i-1$, on d\'eduit $\alpha_j=Q\alpha_j'$. Comme
$\alpha_0',\dots,\alpha_{i-1}'$ sont libres sur $R$, et qu'ils ont les m\^emes degr\'es que ceux d'une base, ils
forment une base de $I_{Y'/X'}$.}

Gruson-Peskine ont montr\'e dans \cite{Peskine}, pour tout degr\'e $\alpha$, parmi les sous-sch\'emas projectifs
 ACM de codimension deux non contenus dans une hypersurface de degr\'e $<d$, l'existence d'un sous-sch\'ema projectif $\Delta$
 {\it minimal}, i.e. v\'erifiant $(\forall i)
\phi_Y(i)\ge
\phi_\Delta(i)$ pour tout autre sous-sch\'ema projectif ACM de codimension deux, non contenu dans une hypersurface de degr\'e
$<d$. On donne ici une variante, lorsqu'on regarde la famille des sous-sch\'emas ACM de codimension deux {\it contenus dans une 
hypersurface irr\'eductible de degr\'e $d$.}¥

\begin{pro} Soit $\alpha=sd-r,r<d$. Soit $\Delta$ le
r\'esiduel d'une intersection compl\`ete $(1,r)¥$¥dans une intersection compl\`ete $(d,s)$. 
Pour tout sous-sch\'ema $Y$ ACM de codimension deux contenu dans une hypersurface irr\'eductible $X$
de degr\'e $d$, on a : $\phi_Y(i)\ge
\phi_\Delta(i)$.
 De plus
$(\forall i)
\phi_Y(i)=
\phi_\Delta(i)$ si et seulement si $Y$ est comme $\Delta$ r\'esiduel d'un
sous-sch\'ema d\'eg\'en\'er\'e de degr\'e $r$ de codimension deux.
\end{pro}

\dem 
Soit $Y$ ACM contenu dans une hypersurface $X$ irr\'eductible de degr\'e $d$.
Alors soit $(n_0,\dots,n_{d-1})$ sa s.c.r. .
Tout d'abord, $n_0\ge s$, sinon $Y$ serait de degr\'e $\le
(s-1)d$ par B\'ezout. 
De plus, pour $i>>0$, $\phi_Y(i)\ge \phi_\Delta(i)$.

 Soit $(n_0',\dots,n_{d-1}')$ la s.c.r. de
$\Delta$ dans une hypersurface irr\'eductible de degr\'e $d$ la contenant. On calcule la s.c.r. de $\Delta$
\`a partir de son r\'esiduel : on enl\`eve, dans la suite $(s,\dots,s+d-1)$, $1$ aux $r$ derniers entiers. 
De $n_{i+1}\le n_{i}+1$, on d\'eduit alors le fait fondamental suivant  : si
$n_i<n_i'$, alors $n_j\le n_j'$ pour tout $j\ge i$. En effet, il n'y a qu'une valeur de $i$ pour laquelle
$n_i'=n_{i+1}'$. 
Donc, le graphe de $(n_{i}¥)$ commence au-dessus de celui de $(n_{i})'$, avant de passer en dessous. 
Consid\'erons la diff\'erence $\Delta^{n-2}\phi_Y(l)- \Delta^{n-2}\phi_\Delta(l)$. On voit que cette diff\'erence, qui est
l'aire s\'eparant les deux graphes en-dessous du niveau $l+1$ (affect\'ee du signe appropri\'e) est d'abord croissante, puis d\'ecroissante,
avant d'\^etre nulle. Elle est donc toujours positive. A fortiori, on en d\'eduit $\phi_Y(i)\ge
\phi_\Delta(i)$. Lorsqu'il y a \'egalit\'e, on voit que $n_{i}=n_{i}'$. En particulier, si on consid\`ere le r\'esiduel $Y'$ de $Y$
dans la section de $X$ par une hypersurface de degr\'e $s$ contenant $Y$, on voit que si $r>0$, $n_{0}(Y')+s+d-2=s+d-1$, i.e. $Y'$ est d\'eg\'en\'er\'e.
\CQFD

On peut par le lemme pr\'ec\'edent red\'emontrer la majoration d'Halphen du genre des courbes gauches:

\begin{pro} Soit $Y$ une courbe alg\'ebrique de degr\'e
 $\alpha$ dans une surface irr\'eductible $X$ de degr\'e $d$
dans
$\P_3$. Ecrivons
$\alpha=sd-r, r<d$. Alors le genre arithm\'etique de $Y$ est inf\'erieur \`a
$G(\alpha,d):=1+sd/2(s+d-4)-r(s+d-r/2-5/2)$. S'il y a \'egalit\'e,  $Y$ est r\'esiduelle d'une courbe plane de degr\'e $r$ dans la
section de $X$ par une surface de degr\'e $s$.
\end{pro}
\dem
Soit ${\bf I}_Y$ le faisceau associ\'e \`a $Y$, et $H$ un plan g\'en\'erique.  La suite exacte: $0\to {\bf I}_Y(l-1)\to
{\bf I}_Y(l)\to {\bf I}_{Y\cap H/H}(l)\to 0$ nous donne :
$$H^0({\bf I}_{Y\cap H/H}(l-1))\to H^1({\bf I}_Y(l-1))\to 
H^1({\bf I}_Y(l))\to H^1({\bf I}_{Y\cap H/H}(l))\to $$
$$H^2({\bf I}_Y(l-1))\to H^2({\bf I}_Y(l))\to 0$$ d'o\`u l'on tire
$h^1({\bf I}_{Y\cap H/H}(l))\le h^1({\bf O}_Y(l-1))-h^1({\bf O}_Y(l))$ puisque $H^2({\bf I}_Y(l))$ est isomorphe \`a $H^1({\bf O}_Y(l))$.
En sommant de $l=1$ \`a $\infty$, on obtient la majoration
$$\sum_l{h^1({\bf I}_{Y\cap H/H}(l))}\le \sum_l{h^1({\bf I}_{\Delta\cap H/H}(l))}=
G(\alpha,d)$$ d'apr\`es la proposition
pr\'ec\'edente. 
De plus, s'il y a \'egalit\'e, on doit avoir que $Y\cap H$ est le r\'esiduel de $r$ points 
align\'es
dans la section de
$C=X\cap H$ par une courbe (dans $H$) de degr\'e $s$. Mais la suite exacte implique alors que : $¥h^1({\bf I}_{Y\cap H/H}(l))=
 h^1({\bf O}_Y(l-1))-h^1({\bf O}_Y(l))$, et donc que :
 
i) $H^1({\bf I}_Y(l))=0$ pour tout $l>0$, puis

ii) $H^0({\bf I}_{Y}(l))\to H^0({\bf I}_{Y\cap H/H}(l))$ est surjectif pour $l\ge 1$. Donc l'existence d'une
courbe de degr\'e $s\ge 1$ contenant
$Y\cap H$ (mais pas $C=X\cap H$) implique l'existence d'une surface $S$ de degr\'e $s$ contenant $Y$ (mais pas $C$, donc pas $X$¥)¥. 
 $Y$ v\'erifie $h^1({\bf I}_{Y¥}¥(l))¥=0$ et est donc ACM.
On voit donc que $Y$ r\'ealise comme $\Delta$ la fonction de Hilbert minimale, et son r\'esiduel dans $X\cap S$ est donc d\'eg\'en\'er\'e.
\CQFD

\section{Syst\`emes lin\'eaires sur les courbes alg\'ebriques et fonctions de Hilbert}

On va s'int\'eresser maintenant plus particuli\`erement au cas o\`u $X$ est de dimension $1$ (i.e. est une courbe
alg\'ebrique) ACM et $Y$ est un groupe de points dessus, particuli\`erement pour l'\'etude des syst\`emes  lin\'eaires sur $X$.
On suppose par la suite pour simplifier que $Y$ est d\'efini localement sur $X$ par une \'equation, i.e. que ${\bf I}_{Y/X}$
 est localement principal.
 Alors, on identifie parfois $Y$ \`a son diviseur de Cartier associ\'e sur $X$, not\'e $[Y]$. Deux groupes de points $Y$ et $Y'$
  de m\^eme degr\'e $\alpha$ sur $X$ sont donc {\it lin\'eairement \'equivalents} si les diviseurs de Cartier associ\'es le sont. 
Il revient au m\^eme de dire qu'il
existe $Z$ sur $X$ tel $Y$ et $Y'$ sont les 
r\'esiduels respectifs de $Z$ dans les sections de $X$ avec des hypersurfaces de m\^eme degr\'e.
 On voit alors en particulier que $Y$ et $Y'$ ont m\^eme
s.c.r. .  Le {\it  syst\`eme lin\'eaire complet} passant pas $Y$ est l'ensemble des groupes de points lin\'eairement \'equivalents \`a $Y$ sur $X$.
On le note $\vline Y\vline$. 
On note ${\bf O}_{X}(Y)$ le faisceau sur $X$ dont la fibre en $x$ est l'ensemble des fonctions rationnelles, 
dont la multiplication par un \'el\'ement de ${\bf I}_{Y,x}$ est
r\'eguli\`ere en $x$.
La {\it dimension} de $\vline Y\vline$ est $h^0({\bf O}_{X}(Y))-1$.

\begin{thm}Soit $X$ est une courbe de Gorenstein. On a $n_{d-1}\le m_{d-1}+s$, avec \'egalit\'e ssi $Y$ est section de $X$ par une 
hypersurface de degr\'e $s$.
\end{thm}
\dem 
Supposons $n_{d-1}\ge m_{d-1}+s$. 
Alors, le terme $(s+m_{d-1}-1-n_{d-1})_{+}¥$ dans $\phi_{Y}(s+m_{d-1}-2)$ est nul, 
et donc $\phi_{Y}(s+m_{d-1}-2)<deg(Y)$, ce qui implique $h^1({\bf I}_{Y}(s+m_{d-1}-2))¥\not= 0$.
D'autre part, $X$ est ACM donc $h^{1}({\bf I}_{X}(l)¥)=0$ pour $l>0$, et donc
 $H^1({\bf I}_{Y}(s+m_{d-1}))\simeq H^1({\bf I}_{Y/X}(s+m_{d-1}))$.
De plus, comme $\omega^1_{X}\simeq {\bf O}_{X}(m_{d-1}-2¥)$, la dualit\'e de Serre nous donne 
$¥H^1({\bf I}_{Y/X}(s+m_{d-1}-2))\simeq H^0({\bf I}_{Y/X}^*(-s))$. 
Mais le faisceau ${\bf I}_{Y/X}^*$ est isomorphe au fibr\'e lin\'eaire d\'efini par le diviseur $[Y]$¥associ\'e \`a $Y$, 
et ${\bf I}_{Y/X}^*(-s)\simeq {\bf O}_X([Y]-[H_s])$, o\`u $H_s$ est une
section de $X$ par une surface de degr\'e $s$. Le diviseur $¥[Y]-[H_s]$¥ est de degr\'e n\'egatif $-r$; il ne
peut avoir de section globale que si il est trivial, auquel cas $Y$ est comme $H_s$ section de
$X$ par une surface de degr\'e $s$, puisque $H^1({\bf I}_{X}(s))=0$. De plus, dans ce cas, on doit avoir $r=0$.
\CQFD

\begin{cor}¥Soit $X$ un sous-sch\'ema projectif irr\'eductible de Gorenstein de dimension $m$ et de degr\'e $d$, de suite caract\'eristique $(m_{i})$.
 Soit $Y\subset X$  une hypersurface ACM, localement principale, de s.c.r. $(n_i)$,  de degr\'e $\alpha:=sd-r, r<d$.
$n_{d-1}\le m_{d-1}+s$, avec \'egalit\'e ssi $Y$ est section de $X$ par une hypersurface de degr\'e $s$.
\end{cor}¥
\dem On se ram\`ene au cas o\`u $X$ est une courbe, en coupant par un sous-espace lin\'eaire g\'en\'erique $H$ de codimension $m-1$.
En effet, les suites caract\'eristiques (absolues et relatives) ne sont alors pas modifi\'ees. En particulier, si $Y\cap H$ a le caract\`ere d'une
section de $X\cap H$ par une hypersurface de degr\'e $s$, alors, $Y$ a la m\^eme s.c.r. dans $X$, et est donc lui-m\^eme section de $X$ par une hypersurface de degr\'e $s$.
\CQFD

Les \'enonc\'es donn\'es ici sur les syst\`emes lin\'eaires reposent sur le lemme fondamental suivant :
\begin{lem} Soit $X$ de Gorenstein. Alors
$dim(\vline Y\vline)=\alpha-\phi_{Y}(m_{d-1}-2)$.
\end{lem}
\dem Le th\'eor\`eme de Riemann-Roch g\'en\'eralis\'e nous donne, pour un diviseur $Y$ sur $X$¥ :
$h^0({\bf O}_{X}(Y))=\alpha+1-p_{a}+i(Y)$, o\`u $p_{a}¥$ est le genre arithm\'etique de $X$ et $i(Y)$ la dimension de
 $H^1({\bf O}_{X}(Y))$, ou encore de $\simeq H^0(\omega^1_{X}(-Y))$ d'apr\`es la dualit\'e de Serre.
Mais ici, comme $X$ est de Gorenstein, $H^0(\omega^1_{X}(-Y))\simeq H^0({\bf O}_{X}(m_{d-1}-2)(-Y))$. Un \'el\'ement de $¥H^0({\bf O}_{X}(m_{d-1}-2))$
appartient \`a $H^0({\bf O}_{X}(m_{d-1}-2)(-Y))$ ssi il s'annule sur $Y$; donc $i(Y)$ \'egale $rg_{k}¥(I_{Y/X}(m_{d-1}-2))$.
On trouve donc bien $dim(\vline Y\vline)=\alpha-\phi_{Y}(m_{d-1}-2)$.
\CQFD

On se donne pour commencer un courbe plane irr\'eductible $X$ de degr\'e $d$, et un groupe de points $Y$ de degr\'e
$\alpha$ sur $X$. On pose
$\alpha:=sd-r$, avec $r<d$. Soit
$\Delta$ le r\'esiduel de $r$ points align\'es dans l'intersection de $X$ avec une courbe de degr\'e $s$. On a vu
que $\phi_\Delta$ est minimale, dans le sens o\`u $\phi_Y(i)\ge \phi_\Delta(i)$ pour tout $i$.
On a par ailleurs toujours  $\phi_Y(i)=\phi_\Delta(i)$ pour $i<s$, ou $i>s+d-3$. La proposition suivante
analyse ce qui se passe dans le cas o\`u l'\'egalit\'e se produit avec $s\le i\le s+d-3$.
 
\begin{pro}
 Si pour un certain entier $i$ compris entre $s$ et $s+d-3$ on a $\phi_Y(i)=\phi_\Delta(i)$, alors :

i) Si $i\ge
s+d-r-1$, alors $n_t=m_t$ pour tout $t\ge d-r$  (alors, $(\forall j\ge
s+d+r-1)\phi_Y(j)=\phi_\Delta(j)$);

ii) Si $i\le s+d-r-3$, on a $n_t=m_t$
pour tout $t\le d-r-1$ (alors, $(\forall j\le s+d-r-1)\phi_Y(j)=\phi_\Delta(j)$);

iii) Si $i=s+d-r-2$, alors:

ou bien
$n_t=m_t$ pour tout $t\le d-r-1$ (auquel cas $(\forall j\le s+d-r-1)\phi_Y(j)=\phi_\Delta(j)$),

ou bien $n_t=m_t$ pour tout $t\ge d-r$ (auquel cas $(\forall j\ge
s+d-r-1)\phi_Y(j)=\phi_\Delta(j)$).
\end{pro}

\dem
Supposons $\phi_Y(i)=\phi_\Delta(i)$ pour un certain $i,s\le i\le s+d-3$. 
On appelle $(n_i)$ la s.c.r. de $Y$ et $(n_i')$ celle de $\Delta$.

Supposons qu'il existe $j$ tel que $n_j\not = n_j'$. On a vu que si $n_j<n_j'$, alors $n_l\le n_l'$
pour $l\ge j$. D'autre part, $\sum_l{(n_l-n_l')}=0$. Le plus petit entier $j$ tel que $n_j\not= n_j'$
doit donc \^etre tel que $n_j>n_j'$. Mais $\sum_l{(n_l-n_l')}=0$ nous montre qu'il existe alors un autre $j$ tel
que $n_j<n_j'$. Si les deux suite $(n_i)$ et $(n_i')$ sont distinctes,
on voit donc que le graphe de $(n_j)$ doit d'abord passer strictement au-dessus de celui de $(n_j')$, puis
strictement en-dessous (mais ne peut plus alors retourner strictement au-dessus). 

D'autre part, si $n_j<n_j'$, on doit m\^eme avoir $n_l<n_l'$ jusqu'\`a $l=d-r$.

Supposons que $i\le s+d-r-2$. Alors,  
de $n_j<n_j'=i$, on d\'eduit $n_{j+1}<n_{j+1}'=i+1$, ce qui implique
$\phi_Y(i+1)<\phi_\Delta(i+1)$, ce qui est impossible d'apr\`es ce qu'on a vu.
Donc : pour $n_j'\le i$, on a $n_j=n_j'$. 

Soit $j_0$ le premier $j$ tel que $n_j\not =n_j'$. On a vu que $n_{j_0}>n_{j_0}'$. 
Mais cela ne peut arriver que si
$j_0\ge d-r$, puisque $n_j$ ne peut avant augmenter plus vite que $n_j'$. On a donc $n_t=n_t'$ pour
$t\le d-r-1$.

On en d\'eduit :
$(\forall j\le s+d-r-2)\phi_Y(j)=\phi_\Delta(j)$.

Supposons maintenant $i\ge s+d-r-1$. De $\phi_Y(i)=\phi_\Delta(i)$, on d\'eduit que la somme
$\sum_j{(i+1-n_j)_+-(i+1-n_j')_+}$ est nulle. Supposons $n_{j_0}<n_{j_0}'$ pour un $j_0\ge d-r$. Alors on aurait
encore
$n_j<n_j'$ pour $j>j_0$. Mais alors la somme pr\'ec\'edente ne peut pas \^etre nulle. On a donc $n_j=n_j'$ pour
$j\ge d-r$.

On en d\'eduit  $(\forall j\ge s+d-r-1)\phi_Y(j)=\phi_\Delta(j)$.

Enfin, supposons $i=s+d-r-2$. 

Le premier cas possible est lorsque $n_j=n_j'$ pour $n_j'\le i+1$. Dans ce cas $n_t=n_t'$ pour $t\le d-r-1$.

Si ce n'est pas le cas, il existe $j_0$ tel que $n_{j_0}'\le i+1$ et $n_{j_0}<n_{j_0}'$. Mais alors $n_j\le n_j'$
pour
$j\ge j_0$â qui implique, comme $\sum_j{(i+1-n_j)_+-(i+1-n_j')_+}$ est nul, que 
$n_{t}=n_{t}'$ pour $t\ge d-r$.

Dans le premier cas, $(\forall j\le s+d-r-2)\phi_Y(j)=\phi_\Delta(j)$, dans le deuxi\`eme $(\forall j\ge
s+d-r-1)\phi_Y(j)=\phi_\Delta(j)$. Les deux cas sont r\'eunis lorsque $\phi_Y=\phi_\Delta$, cas o\`u $Y$ est
comme $\Delta$ r\'esiduel dans la section de $X$ par une courbe de degr\'e $s$ d'un groupe de $r$ points align\'es.
\CQFD

On retrouve par le th\'eor\`eme pr\'ec\'edent la description g\'eom\'erique des syst\`emes lin\'eaires de dimension maximale pour
un degr\'e $\alpha$ donn\'e sur une courbe alg\'ebrique plane de degr\'e $d$,
\'etablie par Ciliberto dans \cite{Ci2} pour les courbes lisses:

\begin{cor}
Soit $\alpha=sd-r$, avec $r<d$. Si $s\ge d-2$ tous les syst\`emes lin\'eaires complets de degr\'e $\alpha$
ont la m\^eme dimension, $\alpha-p$, avec $p=(d-1)(d-2)/2$. Si $s\le d-2$, posons
$r(\alpha)=s(s+3)/2-r$ si $r\le s+1$, et $r(\alpha)=(s-1)(s+2)/2$ si $r\ge s+1$.

Alors la dimension de tout syst\`eme lin\'eaire de degr\'e $\alpha$ est inf\'erieure \`a $r(\alpha)$.
Supposons qu'il passe par $Y$ un syst\`eme lin\'eaire de dimension $r(\alpha)$. Alors :

i) Si $r\le s$, $Y$ est r\'esiduel d'un groupe de $r$ points dans l'intersection de $X$ avec une courbe de degr\'e
$s$;

ii) Si $r\ge s+2$, $Y$ contient l'intersection de $X$ avec une courbe de degr\'e $s-1$;

iii) Si $r=s+1$, 

ou bien $Y$ est r\'esiduel d'un groupe de $r$ points dans l'intersection de $X$ avec une courbe de degr\'e
$s$;

ou bien $Y$ contient l'intersection de $X$ avec une courbe de degr\'e $s-1$.
\end{cor}

\dem
L'\'egalit\'e $dim (\vline Y\vline)=r(\alpha)$ \'equivaut \`a $\phi_Y(d-3)=\phi_\Delta(d-3)$.
Supposons $\alpha\le d(d-3)$.
Si $r\ge s+2$, on d\'eduit de la proposition pr\'ec\'edente que $n_{d-1}=s+d-2,n_{d-2}=s+d-3,\dots,
n_{d-r}=s+d-r-1,\dots$. Cela signifie, comme $s+d-r\le d-2$, que la suite caract\'eristique absolue
de $Y$, obtenue en supprimant les doubles, a un trou entre $n_{d-s+1}'=n_{d-s+1}=d$ et $n_{d-s}\le d-2$.

D'apr\`es ce qui pr\'ec\`ede, on en d\'eduit que $Y$ contient un groupe de points $Y'$, section de $Y$ avec une
courbe $X'$ de degr\'e $s-1$, $Y'$ ayant pour s.c.r. $(s-1,\dots)$. Mais la section de $X$ par $X'$
contient $Y'$ et a la m\^eme s.c.r. , donc est \'egale \`a $Y'$. Donc $Y$ contient $Y'$, section de
$X$ avec une courbe $X'$ de degr\'e $s-1$.

Si $r\le s$, on d\'eduit de la proposition pr\'ec\'edente que $n_j=n_j'$ pour $j\le d-r-1$. Donc $n_0=s$, et $Y$
est contenu dans une courbe de degr\'e $s$, donc dans la section de cette courbe de degr\'e $s$ avec $X$. 

Si $r=s+1$, il y a deux cas possibles.

Dans le premier cas, $n_j=n_j'$ pour $j\le d-r-1$. Alors $Y$ est contenu dans la section de $X$ avec une courbe de
degr\'e
$s$.

Dans le deuxi\`eme cas,
$n_j=m_j$ pour
$j\ge d-r+1$. Dans ce deuxi\`eme cas, on voit pour la m\^eme raison que ci-dessus que $Y$ contient l'intersection
de $X$ avec une courbe de degr\'e $s-1$.
\CQFD

{{\bf Remarque}.
1. Lorsque $r=s+1$, $Y$ satisfait les deux conditions: contenir la section de $X$ avec une courbe de degr\'e $s-1$,
et
\^etre contenu dans la section de $X$ avec une courbe de degr\'e $s$, ssi $Y$ est comme $\Delta$, r\'esiduel de $r$ points align\'es dans
la section de $X$ avec une courbe de degr\'e $s$. 

2. Supposons $Y$ contient la section de $X$ avec une courbe de degr\'e $s-1$, et soit $Y'$ le r\'esiduel (de
degr\'e $d-r$) de cette section dans $Y$. Alors, $Y'$ est la partie fixe du syst\`eme lin\'eaire
$\vline Y \vline$. 
Supposons $Y$ contenu dans la section de $X$ par une courbe de degr\'e $s$, et soit $Y''$ le r\'esiduel (de
degr\'e $r$) de $Y$ dans cette section.  Si
$r\le s$, le syst\`eme lin\'eaire 
$\vline Y\vline$ n'a pas de point fixe.  Si $r=s+1$ mais que $Y''$ n'est pas align\'e,  alors $\vline Y\vline$
n'a pas de point fixe.  }

Soit $X\subset \P_n$ une courbe alg\'ebrique irr\'eductible. On voudrait voir ce qui arrive, lorsqu'on "ajoute" \`a
un groupe de points
$Y'\subset X$ un autre groupe de points
$Y''\subset X$, pour obtenir un groupe de points $Y\subset X$, comment se transforme la s.c.r. lorsqu'on passe de
$Y'$ ($(n_i')$) \`a $Y$ ($(n_i)$), en fonction de $Y''$.

Observons qu'on peut ajouter \`a $Y'$ respectivement deux groupes $Y''_1$ et $Y''_2$, avec des s.c.r.
distinctes, mais obtenir pour
$Y_1=Y'\cup Y''_1$ et $Y_2=Y'\cup Y''_2$ les m\^emes s.c.r. . On ne peut donc pas en
g\'en\'eral calculer la s.c.r. de $Y''$ \`a partir de celles de $Y=Y'\cup Y''$ et $Y'$, comme on l'a
fait lorsque $Y'$ est une section de $X$ par une hypersurface de $\P_n$.

En g\'en\'eral, lorsque $Y''$ est un point, i.e. qu'on passe de $Y'$ \`a un groupe de points
$Y$ sur $X$ contenant $Y'$ dont le degr\'e est plus grand d'une unit\'e, on "ajoute une case" sur le graphe de la fonction
$i\mapsto n_i'$, graphe que l'on peut voir comme une superposition de cases, \`a un certain niveau (si plusieurs
valeurs successives de la suite $(n_i)$ sont \'egales \`a un entier $l$, l'ajout de case au niveau $l$ se fera pour
la derni\`ere valeur de $i$ pour laquelle $n_i=l$). Cet ajout doit se faire de sorte que l'in\'egalit\'e
$n_{i+1}\le n_i+1$ reste v\'erifi\'ee; seuls certains "ajouts de case" correspondent \`a un "ajout de point".

On d\'efinit pour chaque degr\'e $i$ le groupe de points $Y_i$ sur $X$ d\'efini par l'id\'eal de $A_X$, contenu dans
$I_{Y/X}$, que l'on obtient en ne conservant comme g\'en\'erateurs que les polyn\^omes homog\`enes de
$I_{Y/X}$ de degr\'e $\le i$. On a donc $Y_{n_0+s}=Y\subset\dots\subset Y_{n_0+1}\subset
Y_{n_0}$, o\`u l'on suppose que
$I_{Y/X}$ est engendr\'e par des polyn\^omes de degr\'e $\le n_0+s$.

\begin{pro}
Un ajout de case sur le niveau $n_0+i+1$ correspond \`a "ajouter un point \`a $Y$" ssi
$Y_{n_0+i+1}\subset Y_{n_0+i}$ est une inclusion stricte; on peut alors ajouter une case au niveau
$n_0+i+1$ en ajoutant \`a $Y$ un point de $Y_{n_0+i}-Y_{n_0+i+1}$.
En particulier, si l'on ajoute \`a $Y$ un point en dehors de $Y_{n_0}$, on ajoute une case sur le niveau
de base $n_0$.
\end{pro}

\dem
Soit $\alpha_0,\alpha_1,\dots,\alpha_{d-1}$ les g\'en\'erateurs de $I_{Y/X}$ comme $R_{1}¥-$module.
Alors $I_{Y_{n_0+i}/X}(l)=I_{Y/X}(l)$ pour $l\le n_0+i$.
Supposons que $Y_{n_0+i+1}\subset Y_{n_0+i}$ est une inclusion stricte.
Soit $Y'$ un groupe de points de degr\'e $deg(Y)+1$, contenu dans
$Y_{n_0+i}$ mais pas dans $Y_{n_0+i+1}$. 

Il existe un polyn\^ome de degr\'e $n_0+i+1$, s'annulant sur $Y_{n_0+i+1}$ mais pas sur $Y'$. D'autre part, en
degr\'e
$l\le n_0+i$, les polyn\^omes de $I_{Y/X}$, $I_{Y_{n_0+i}/X}$, et $I_{Y'/X}$, sont les m\^emes. 

Soit
$(n_i')$ la s.c.r. de $Y'$. On a :
$rg_{k} (I_{Y'/X} (l)) = rg_{k}(I_{Y/X}(l))$,
donc $\sum_{i=0}^{d-1}{(l+1-n_i)_+}=\sum_{i=0}^{d-1}{(l+1-n_i')_+}$ pour $l=n_0+i$.
D'autre part, pour $l=n_0+i+1$, l'inclusion $I_{Y'/X}(l)\subset I_{Y/X}(l)$ est stricte, donc
$$\sum_{i=0}^{d-1}{(l+1-n_i)_+}>\sum_{i=0}^{d-1}{(l+1-n_i')_+}.$$
Donc, on passe de $(n_i)$ \`a 
$(n_i')$ en ajoutant une case sur le niveau $n_0+i+1$.

D'autre part, supposons qu'on ait l'\'egalit\'e $Y_{n_0+i+1}= Y_{n_0+i}$. 
Soit $Y'$ un groupe de points obtenu \`a partir de $Y$ en lui ajoutant un point, et tel que $Y'_{n_{0}+i}=Y_{n_{0}+i}$.
Alors $Y_{n_{0}+i+1}\subset Y'_{{n_{0}¥+i+1}¥}¥\subset Y'_{n_{0}+i}=Y_{n_{0}+i}=Y_{n_{0}+i+1}$.
On a donc $Y_{n_{0}+i+1}'=Y_{n_{0}+i+1}$, et donc on ne peut pas ajouter de case au niveau $n_{0}+i+1$.
\CQFD

Pour ajouter une case sur le niveau $n_j$, il faut d'apr\`es l'in\'egalit\'e $n_{i+1}\le n_i+1$, si $j>0$, que $l_{n_{j}¥}\ge 2$. 
 Cette condition n\'ecessaire n'est pas toujours suffisante. N\'eanmoins on peut montrer:

\begin{lem} Soit $(n_i)$ la s.c.r. de $Y$. Soit $j$ le premier entier tel que
$n_j=n_{j+1}$. Alors il existe un groupe de point $Y'$ sur $X$ contenant $Y$, de degr\'e $deg(Y)+1$, tel que la
s.c.r. $(n_i')$ de $Y'$ soit obtenue \`a partir de celle de $Y$ en ajoutant une case sur le niveau
$n_j$.
\end{lem}

\dem
Soit $j$ l'entier donn\'e dans l'\'enonc\'e. Alors on peut choisir $I_{Y/X}$,
$\alpha_0$ (de degr\'e $n_0$), $\alpha_1=Y_2\alpha_0$,\dots,$\alpha_{j-1}=Y_2^{j-1}\alpha_0$. Ainsi, en degr\'e
$n_{j-1}$, l'id\'eal $I_{Y/X}$ est engendr\'e par $\alpha_0$. En degr\'e $n_j$, ce n'est plus le cas  puisque
$n_j=n_{j+1}$. Ainsi, $Y_{n_j}$ est strictement inclus dans $Y_{n_0}$ (tout en contenant $Y$). On peut donc trouver
un groupe de points $Y'$ de degr\'e $deg(Y)+1$, contenant $Y$, contenu dans $Y_{n_0}$, mais pas contenu dans
$Y_{n_j}$. Alors, la s.c.r. de $Y'$ est obtenue \`a partir de celle de $Y$ en ajoutant une case
sur le niveau $n_j$.
\CQFD

{{\bf Remarque.} Ce n'est pas parcequ'un niveau est de largeur $>1$ qu'on peut toujours rajouter une case dessus.
Soit $X$ une sextique plane. Soit $Y$ donn\'e par $5$ points align\'es sur $X$ et $4$ points g\'en\'eriques sur
$X$. $Y$ a alors comme s.c.r.
$(3,3,4,4,5,5)$. La s.c.r. $(3,3,4,5,5,5)$ est obtenue par la r\'eunion $Y'$ de $9$ points
sur une section conique de $X$ et d'un point g\'en\'erique de $X$. Il n'est pas possible que $Y$ soit contenu dans
$Y'$.
Donc, la possibilit\'e d'ajouter un point sur un niveau (qui, s'il n'est pas le "niveau de base" $n_0$, doit
\^etre de largeur $\ge 2$ pour cette possibilit\'e) nous donne de l'information sur
$Y$. Dans l'exemple pr\'ec\'edent, o\`u $X$ est une sextique plane, si $Y$ est la r\'eunion de $8$ points d'une
section conique et d'un point g\'en\'erique de $X$, $Y$ a encore comme s.c.r.
$(3,3,4,4,5,5)$ et il est possible de former la s.c.r. 
$(3,3,4,5,5,5)$ en ajoutant \`a $Y$ l'un des $4$ points restants de la section conique.
}

Soit $(n_i)$ la s.c.r. de $Y$ dans $X$. On a vu que $n_i\ge i$, et si $X$ est
irr\'eductible, $n_{i+1}\le n_i+1$. On peut montrer le th\'eor\`eme suivant :

\begin{thm} Soit $X$ une courbe irr\'eductible de $\P_2$, de degr\'e $d$. Pour toute suite $(n_i)_{0\le
i\le d-1}$ v\'erifiant
$n_i\ge i, n_i\le n_{i+1}\le n_i+1$, on peut construire sur $X$ un groupe de points $Y$ sur $X$ ayant $(n_i)$ pour
s.c.r. .
\end{thm}
\dem
La d\'emonstration se fait par r\'ecurrence sur la somme
$\sum_{i=0}^{d-1}{(n_i-i)}$, le {\it degr\'e} de la suite $(n_{i}¥)¥$¥. L'ensemble vide r\'ealise la suite $n_i:=i$. Supposons que toutes les suites
v\'erifiant $n_i\ge i, n_i\le n_{i+1}\le n_i+1$,
$\sum_{i=0}^{d-1}{n_i-i}=\alpha$ soient r\'ealis\'ees. Soit une suite $(n_i)$ de degr\'e $\alpha+1$. On consid\`ere
le premier entier $i$ tel que $n_{i+1}=n_i$. S'il n'y en a pas, la suite est de la forme $n_i=n_0+i$; elle est
r\'ealis\'ee par la section de $X$ avec une courbe de degr\'e $n_0$ la coupant proprement. Sinon, un tel entier
existe, on l'appelle $j$; on construit une nouvelle suite $n_i'$ en posant $n_i'=n_i$, sauf si $i=j$, o\`u
$n_j':=n_j-1$.

 Alors $\sum_i{(n_i'-i)}=\alpha$, et d'apr\`es l'hypoth\`ese de r\'ecurrence on peut
r\'ealiser la suite $(n_i')$ pour un groupe de points $Y'\subset X$. D'apr\`es le lemme pr\'ec\'edent, on peut en
ajoutant \`a $Y'$ un point, obtenir un groupe de points $Y$ dont la s.c.r. est pr\'ecis\'ement $(n_i)$, ce qui
termine la d\'emonstration.
\CQFD

\section{Questions}

Soit $X$ une courbe plane de degr\'e $d$ $m_{d-1}=d-1$, et $Y$un groupe de points localement principal sur $X$.
Pour tout $s$, soit $r(sd)$ la dimension du syst\`eme lin\'eaire de degr\'e $sd$ d\'efini par les sections de $X$ avec les courbes de degr\'e $s$. On a
$r(s)=s(s+3)/2$ si $s<d$. 
On a montr\'e ci-dessus l'\'enonc\'e pr\'ec\'edent:

{\it Pour tout $s\le m_{d-1}-2$, et tout $Y$ de degr\'e $sd$ sur $X$,  $dim(\vline Y\vline)\le r(sd)$, avec \'egalit\'e ssi $Y$ est la section
de $X$ avec une hypersurface de degr\'e $s$.
}

On demande si cet \'enonc\'e reste valide dans le cadre plus g\'en\'eral o\`u $X$ est une courbe de Gorenstein dans $\P_{n}$. Il le serait si on pouvait d\'emontrer
la conjecture suivante:

{{\it Conjecture.}

Soit $X\subset \P_{n}$ de Gorenstein. 
Soit $\Delta$ la section de $X$ par une hypersurface de degr\'e $s$. Alors, pour tout groupe de points localement principal $Y$ de degr\'e $sd$,
on a $\phi_{Y}(l)\ge \phi_{\Delta}(l)¥$ pour tout $l$. De plus, l'ensemble des $l$ tels que $\phi_{Y}(l)\not=\phi_{\Delta}(l)$ est connexe.
}

Il d\'ecoule de ce qui pr\'ec\`ede que si $\phi_{Y}(l)= \phi_{\Delta}(l)¥$ pour tout $l$, alors la s.c.r. $(n_{i})$ de $Y$ est la m\^eme que celle 
de $X$ et donc $Y$ est la section de $X$ par une hypersurface de degr\'e $s$. 

Supposons que la conjecture est v\'erifi\'ee.
Soit donc $Y$ de degr\'e $sd$ sur $X$, avec $s\le m_{d-1}-2$. Si $dim(\vline Y\vline)=r(sd)$, alors $\phi_{Y}(m_{d-1}-2)=\phi_{\Delta}(m_{d-1}¥-2)$
d'apr\`es ce qu'on a vu. Mais alors, les fonctions de Hilbert sont \'egales avant ou apr\`es,  puisque l'ensemble des $l$ 
tels que $\phi_{Y}(l)\not=\phi_{\Delta}(l)$ est connexe. Si elles le sont avant, $Y$ v\'erifie $n_{0}=s$, et donc $Y$ est section de $X$ par une hypersurface de
degr\'e $s$. Si elles le sont apr\`es, alors, $n_{d-1}=m_{d-1}+s$. Mais alors, on a vu ci-dessus que dans ce cas aussi, $Y$ est section de $X$ par une hypersurface de
degr\'e $s$.

On esp\`ere pouvoir trouver des propri\'et\'es de la s.c.r. $n_{i}$ qui nous permettent de montrer la conjecture pr\'ec\'edente.

\appendix{\bf Appendice}

{ 
{\bf Sous-sch\'emas li\'es et r\'esiduel}
 
 Soit $X$ et $X'$ deux c\^ones de $A^{n+1}$; on suppose $Z=I_X\cap I_{X'}$. Alors on a un morphisme 
 de $A_n-$modules 
$I_X'\to Hom_{A_n}(I_X,I_Z)$. On dit que $X$ et $X'$ sont li\'es, et que $X'$ est
 {\it r\'esiduel} de $X$ dans $Z$, si ce morphisme est un isomorphisme.

{\bf Th\'eor\`eme des syzygies gradu\'e}

 Etant donn\'e un $A_n-$module gradu\'e de type fini $M$, on consid\`ere une suite exacte 
$$0\to L\to
\oplus_{i_s=0}^{d_s-1}{A_n[-i_s]}\to \dots \to \oplus_{i_0=0}^{d_0-1}{A_n[-i_0]}\to M\to 0$$
 Alors, le th\'eor\`eme des syzygies gradu\'e dit que si $s\ge n$,
$L$ est un $A_n-$module libre gradu\'e.

On en d\'eduit:
\begin{lem}
Soit $X$ un sous-sch\'ema projectif ACM de dimension $m$. On se donne $m+r$ formes lin\'eaires lin\'eairement ind\'ependantes
$Y_0,\dots,Y_{m+r}$, dont l'annulation d\'efinit un sous-espace
projectif ne rencontrant pas $X$. Soit
$R_{m+r}=k[Y_0,\dots,Y_{m+r}]$. Alors pour toute suite exacte :
$$0\to L\to \oplus_{i=0}^{d_s-1}{R_{m+r}[-m_{i,s}]}\to \dots \to \oplus_{i=0}^{d_0-1}{R_{m+r}[-m_{i,0}]}
\to A_X\to 0$$ avec
$s\ge r$, $L$ est un $R_{m+r}$-module libre.
\end{lem}

\dem 
On fait r\'ecurrence sur la dimension $m$. 
Le cas de $m=-1$ d\'ecoule directement du th\'eor\`eme des syzygies gradu\'e. 
Supposons que le lemme soit vrai en
dimension $m-1$. Soit
$Z_0,\dots, Z_m$ $m+1$ combinaisons $k-$lin\'eaires des $Y_i$ ne rencontrant pas le support de $X$ (il suffit de choisir des 
combinaisons $k-$lin\'eaires g\'en\'eriques).
Alors, $(Z_1,\dots,Z_m)$ est une suite r\'eguli\`ere dans
$A_X$. On applique
\`a la suite exacte $$0\to L\to \oplus_{i_s=0}^{d_s-1}{R_{m+r}[-i_s]}\to \dots \to \oplus_{i_0=0}^{d_0-1}{R_{m+r}[-i_0]}
\to A_X\to 0$$ le foncteur
$\otimes_{R_{m+r}}{R_{m+r}/Z_m R_{m+r}}$; la multiplication par $Z_m$ est
injective dans $A_X$. Alors, un calcul de rang sur $k$ des $k-$espaces vectoriels nous montre que bien
que le foncteur ne soit pas exact \`a gauche, il conserve ici la suite exacte. De plus, $(Z_0,\dots,Z_{m-1})$
 forme une suite r\'eguli\`ere pour $X\cap \{Z_m=0\}$. On peut donc appliquer l'hypoth\`ese de r\'ecurrence sur $m$.
  Le fait que $L/Z_{m}L$ soit libre sur $R_{m+r}/Z_{m}R_{m+r}$ nous montre alors que $L$ est libre sur $R_{m+r}$.
\CQFD

 }

\end{document}